\documentclass[12pt]{amsart}
\usepackage{amsmath}
\usepackage{amsthm}
\usepackage{amssymb}
\usepackage{mathrsfs}
\usepackage{amscd,verbatim}
\usepackage[all]{xy}
\usepackage{bm}
\usepackage{url}
\usepackage{enumitem}
\usepackage{hyperref}
\usepackage[backend=biber,style=alphabetic,doi=false,isbn=false,url=false,eprint=false]{biblatex}
\theoremstyle{definition}
\newtheorem{thm}{Theorem}[section]
\newtheorem{definition}[thm]{Definition}
\newtheorem{prop}[thm]{Proposition}
\newtheorem{lem}[thm]{Lemma}
\newtheorem{rem}[thm]{Remark}
\newtheorem{cor}[thm]{Corollary}

\newtheorem*{ack}{Acknowledgements}
\numberwithin{equation}{section}
\usepackage[top=20truemm, bottom=20truemm, left=20truemm, right=20truemm]{geometry}

\newcommand{\Z}{\mathbb{Z}}
\newcommand{\Q}{\mathbb{Q}}

\newcommand{\Spec}{\operatorname{Spec}}

\newcommand{\Hom}{\operatorname{Hom}}

\newcommand{\Ext}{\operatorname{Ext}}

\newcommand{\Frac}{\operatorname{Frac}}

\newcommand{\End}{{\rm End}}

\newcommand{\Ker}{\operatorname{Ker}}
\newcommand{\im}{\operatorname{Im}}
\newcommand{\Cok}{\operatorname{Coker}}

\newcommand{\uni}{{\rm uni}}
\newcommand{\red}{\mathrm{red}}
\newcommand{\MT}{\mathrm{m}}
\newcommand{\LL}{\mathrm{ll}}

\newcommand{\scrO}{\mathscr{O}}

\newcommand{\calC}{\mathcal{C}}

\newcommand{\calF}{\mathcal{F}}

\newcommand{\bbP}{\mathbb{P}}

\newcommand{\Def}{:=}
\newcommand{\et}{\text{{\rm \'et}}}

\newcommand{\ab}{\mathrm{ab}}

\newcommand{\fppf}{{\rm fppf}}

\newcommand{\G}{\mathbb{G}}
\newcommand{\F}{\mathbb{F}}
\newcommand{\A}{\mathbb{A}}

\newcommand{\Ab}{\mathrm{Ab}}

\newcommand{\frkm}{\mathfrak{m}}

\newcommand{\Nori}{\mathrm{N}}
\newcommand{\Pro}{\operatorname{Pro}}

\newcommand{\Ind}{\operatorname{Ind}}
\newcommand{\Diag}{\operatorname{Diag}}
\newcommand{\Pic}{\mathrm{Pic}}

\newcommand{\cond}{\mathfrak{c}}
\newcommand{\fil}{\mathrm{fil}}
\newcommand{\cl}{\mathrm{cl}}
\newcommand{\pr}{\mathrm{pr}}

\usepackage[usenames,dvipsnames]{color}

\usepackage[T1]{fontenc}

\usepackage[usenames,dvipsnames]{color}

\begin{document}

\title{Class field theory for function fields and finite abelian torsors}

\author[B. Cais]{Bryden Cais}
\address{Department of Mathematics, The University of Arizona,
Tucson, Arizona, USA}

\email{cais@arizona.edu}

\author[S. Otabe]{Shusuke Otabe}
\address{Graduate School of Engineering, Nagoya Institute of Technology, Nagoya, Aichi, Japan
}
\email{shusuke.otabe@nitech.ac.jp}


\date{\today}

\keywords{geometric class field theory, fundamental group schemes, curves, positive characteristic}
\subjclass{14H30, 14H40, 14L15}

\thispagestyle{empty}

\begin{abstract}
Let $U$ be a smooth and connected curve over an algebraically
closed field of positive characteristic, with smooth compactification $X$.
We generalize classical Geometric Class Field theory
to provide a classification of fppf $G$-torsors 
over $U$ in terms of isogenies of generalized Jacobians, for any finite
abelian group scheme $G$.  
We then apply this classification to give a novel description 
of
the abelianized Nori fundamental group scheme 
of $U$ in terms of the Serre--Oort fundamental groups
of generalized Jacobians of $X$; when $U=X$ is projective, 
we recover a well-known description of the abelianized fundamental group scheme of $X$ as the projective limit of all torsion subgroup schemes of
its Jacobian.
\end{abstract}

\maketitle

\thispagestyle{empty}

\section{Introduction}\label{sec:intro}

Pioneered by Lang \cite{Lang} and Rosenlicht \cite{Rosenlicht54,Rosenlicht57} in the 1950s,
and 
expounded in Serre's book \cite{Serre59},
Geometric Class Field theory provides a beautiful description of the 
abelian covers of 
a smooth, projective 
and connected curve $X$ over an algebraically closed field $k$ 
in terms of \'etale isogenies onto its generalized Jacobians.
In the case of unramified coverings, this description is encapsulated
by an identification of profinite abelian groups
\begin{equation}
    \pi_1^{\et}(X)^{\ab} \simeq \varprojlim_n J_X[n](k),
    \label{eq:Lang}
\end{equation}
with $J_X$ the Jacobian of $X$.  Via Nori's theory of the fundamental group scheme \cite{Nori82},
\eqref{eq:Lang} admits a striking enrichment to an isomorphism
of profinite group {\em schemes}
\begin{equation}
    \pi^{\Nori}(X)^{\ab} \simeq \varprojlim_n J_X[n]
    \label{eq:Nori}
\end{equation}
recovering \eqref{eq:Lang} on $k$-points; see
\cite[Remark 2]{Nori83} and cf.~\cite[Corollary 3.8]{Antei}, \cite[Corollary 7.2]{Langer}.  The isomorphism \eqref{eq:Nori} is likewise equivalent to
an enrichment of unramified Geometric Class Field theory
to a classification of {fppf} $G$-torsors over $X$
for an arbitrary finite abelian $k$-group scheme $G$
in terms of (not necessarily separable) isogenies onto $J_X$.

It is natural to ask for analogues of \eqref{eq:Lang}--\eqref{eq:Nori}
in the case of {\em ramified} coverings, or equivalently, 
for analogous descriptions of the abelianized fundamental group
(scheme) of a smooth {\em affine} curve.  
In the case of
\eqref{eq:Lang}, such a description essentially goes back to Serre \cite[\S5.2]{Serre61},
whose work gives an identification of profinite abelian groups
\begin{equation}
        \pi_1^{\et}(X\setminus S)^{\ab} \simeq \varprojlim_{\frkm} \varpi_1(J_{X,\frkm})(k)
        \label{eq:serre-affine}
\end{equation}
for any finite set $S\subseteq X(k)$, 
in which $\varpi_1(\cdot)$ is the {\em Serre--Oort fundamental group}, as defined by
Oort \cite[II,\S7]{Oort} (cf.~\cite{Brion}), $J_{X,\frkm}$
is the generalized Jacobian of $X$ with modulus $\frkm$,
and
the inverse limit is taken over all $\frkm$ supported on $S$.

In this paper, we establish the natural analogue of \eqref{eq:Nori}
for arbitrary smooth affine curves by enriching \eqref{eq:serre-affine} 
to an isomorphism of profinite group schemes:

\begin{thm}\label{thm:main-intro}
    Let $X$ be a smooth projective and connected curve over $k$
    and $S\subseteq X(k)$ any finite set of points.  There is an isomorphism
    of profinite group schemes
    \begin{equation}
        \pi^{\Nori}(X\setminus S)^{\ab} \simeq \varprojlim_{\frkm} \varpi_1(J_{X,\frkm}),
        \label{eq:nori-affine}
\end{equation}
    recovering \eqref{eq:serre-affine} on $k$-points.
\end{thm}

\begin{rem}
    A precursor to the Serre--Oort fundamental group $\varpi_1(\cdot)$
    was defined by Serre \cite[\S5.3]{Serre_GP}, who constructed
    a functor from 
    the category of
    proalgebraic $k$-groups (in the sense of \cite[\S2.1 Definition I]{Serre_GP}) to the category of profinite abelian groups.  At the time, Grothendieck's theory of group {\em schemes} was in its infancy, so Serre's category of proalgebraic groups may look strange
    to a modern reader: it is the pro-category of the category of {\em quasi-algebraic groups},
    the latter being the localization of the category of commutative, smooth algebraic $k$-groups
    at the class of morphisms inducing an isomorphism on $k$-points. 
    Oort \cite[II.7]{Oort} generalized Serre's construction
    to define $\varpi_1(\cdot)$ as a functor from the pro-category of commutative algebraic $k$-groups
    to the category of profinite commutative $k$-group {\em schemes}.  We review Oort's construction 
    in \S\ref{sec:ab Nori}, noting here that for a commutative and connected algebraic $k$-group $H$, one has the ``concrete'' description
    \begin{equation}
            \varpi_1(H) = \varprojlim_{\varphi: H'\twoheadrightarrow H} \Ker(\varphi)
            \label{eqn:varpi-concrete}
    \end{equation}
    with the limit taken over all isogenies  $\varphi:H'\twoheadrightarrow H$
    of commutative algebraic $k$-groups.  Thus, \eqref{eq:nori-affine} provides a description of the abelianized
    Nori fundamental group of $X\setminus S$ in terms of isogenies onto the generalized Jacobians of 
    $X$ with moduli supported on $S$.
    More or less by definition, the associated profinite abelian
    group of $k$-points $\varpi_1(\cdot)(k)$ recovers Serre's original construction. 

    Note that when $S=\emptyset$,
    \eqref{eq:serre-affine} and \eqref{eq:nori-affine} recover
    \eqref{eq:Lang} and \eqref{eq:Nori}, respectively, 
    as every isogeny with target $J_X$ factors multiplication
    by $n$ for some $n$, by duality for (isogenies of) 
    abelian varieties.
\end{rem}

Let us briefly describe the proof of Theorem \ref{thm:main-intro}.
Writing $U:=X\setminus S$, Nori's work provides a profinite commutative $k$-group scheme 
$\pi^{\Nori}(U)^{\ab}$ with a functorial identification
\begin{equation}
    \Hom(\pi^{\Nori}(U)^{\ab},G) = H^1_{\fppf}(U,G)
\end{equation}
for variable finite commutative $k$-group schemes $G$.  On the other hand, 
it follows from the concrete description of $\varpi_1(\cdot)$ given by \eqref{eqn:varpi-concrete}
that, for any commutative, smooth and connected $k$-group $H$,
pushout of the tautological extension 
arising from an isogeny $\varphi:H'\rightarrow H$
along a homomorphism of finite $k$-group schemes $\Ker(\varphi)\rightarrow G$
induces a functorial identification (see Lemma \ref{lem:univ cov B})
\begin{equation}
    \Hom(\varpi_1(H),G) = \Ext^1(H,G). 
\end{equation}
As such, the identification \eqref{eq:nori-affine} is equivalent to a functorial identification
\begin{equation}
    H^1_{\fppf}(U,G) \simeq \varinjlim_{\frkm} \Ext^1(J_{X,\frkm},G).
    \label{eqn:fppf-ext}
\end{equation}
In the case $S=\emptyset$, so $U=X$ is projective, \eqref{eqn:fppf-ext}
follows from the well-known identifications
\begin{equation}\label{eq:Raynaud}
H^1_{\fppf}(X,G)\simeq\Hom(G^D,\underline{\Pic}_{X})\simeq\Hom(G^D,\underline{\Pic}_{J_X})\simeq H^1_{\fppf}(J_X,G)\simeq\Ext^1(J_X,G),
\end{equation}
where $G^D$ is the Cartier dual of $G$; see \cite[Proposition 6.2.1]{Raynaud}, \cite[III, \S4, Proposition 4.16]{Milne} for the first and third isomorphisms; the second isomorphism is due to the autoduality of the Jacobian $J_X$; for the last isomorphism, see \cite[III, \S4, Corollary 4.20]{Milne}. 

Thus, the main contribution of this paper is to establish \eqref{eqn:fppf-ext}
for general $S\subseteq X(k)$.  We will prove a more precise result: fixing 
$x_0\in U(k)$ and, for any modulus $\frkm$ supported on $S$, writing $\varphi_{\frkm}:U\rightarrow J_{X,\frkm}$
for the corresponding Albanese map, we establish the following generalization of Geometric Class Field theory:

\begin{thm}(cf.\ Theorem \ref{thm:CFT})\label{thm-intro:CFT}
    For any finite abelian $k$-group scheme $G$ and any {fppf} $G$-torsor
    $P\rightarrow U$, there exists a modulus $\frkm$ supported on $S$ 
    and a unique isogeny 
    of algebraic groups $J'\rightarrow J_{X,\frkm}$ with kernel $G$ such that $P=\varphi_{\frkm}^* J'$ is the pullback
    of $J'$ along $\varphi_{\frkm}$.
\end{thm}

When $G$ is reduced, Theorem \ref{thm-intro:CFT} recovers the classical Geometric Class Field theory
of Lang, Rosenlicht, and Serre; see \cite{Serre}.  In fact, the identification 
\eqref{eqn:fppf-ext} is stated by Deligne for arbitrary commutative $k$-groups $G$
in a 1974 letter to Serre (reproduced in the Appendix of \cite{BE}).
Deligne sketches a proof, which proceeds by d\'evissage to the case of smooth $G$.
He then interprets each side of \eqref{eqn:fppf-ext} as the $i=1$
degree of a cohomological $\delta$-functor in {\em two} variables, and argues that each $\delta$-functor
is effaceable and therefore universal by (a generalization of) Grothendieck's theorem \cite[II, Proposition 2.2.1]{Grothendieck}.  In this way, Deligne's approach reduces \eqref{eqn:fppf-ext} to 
the identification $\Hom((U,x_0),(G,0))\simeq \varinjlim_{\frkm} \Hom(J_{X,\frkm},G)$
for smooth commutative $k$-groups $G$, where $\Hom((U,x_0),(G,0))$ stands for the set of pointed maps, which follows from Rosenlicht's work
(see \cite[V, \S2, No.\,9, Theorem 2]{Serre}).  While beautiful, Deligne's sketch 
is missing some details (e.g.~we are unaware of any reference for the
generalization of Grothendieck's theorem to the case of cohomological $\delta$-functors in
two variables that the argument relies on).  Rather than fill in these details, we provide
a different---and in some ways more elementary---approach, establishing Theorem \ref{thm-intro:CFT}
in \S\ref{sec:abelian torsors}.

Among all moduli $\frkm$
guaranteed by Theorem \ref{thm-intro:CFT}, there is a unique {\em minimal} one (see Corollary  \ref{cor:minmod}). 
In the classical case of \'etale (i.e.~constant) $G$,
this minimal modulus coincides with the conductor divisor of the associated branched covering of curves.
In the general case, it yields a natural notion
of $G$-torsor over $U$ {\em with ramification no worse than $\frkm$}.
Such a notion has recently been proposed by R\"ulling and Saito 
\cite[\S9]{RS} for smooth schemes over a perfect field as an application of their general theory of ramification for reciprocity sheaves; see also Remark \ref{rem:RS}. We will give a description of the minimal modulus purely in terms of local fields. To that end, we introduce a conductor $\cond^G_{\Frac\widehat{\scrO}_{X,x}}(-)$ (see Definition \ref{def:cond}) for finite local-local torsors over the fraction field of a completed local ring. 
In \S\ref{sec:torsors with modulus}, we provide the following description of the minimal modulus $\frkm$ for finite local-local torsors.

\begin{thm}(cf.\,Theorem \ref{thm:m vs c})\label{thm-intro:m vs c}
Let $G$ be a finite local-local $k$-group and $P\in H^1_{\fppf}(U,G)$. Let  $\frkm(P)=\sum_{x\in S}n_x(P) x$ be the minimal modulus of $P$. Then for any $x\in S$, we have
\[
n_x(P)=\cond^G_{\Frac\widehat{\scrO}_{X,x}}(P|_{\Frac\widehat{\scrO}_{X,x}}).
\]
\end{thm}

We also study the particular case of $G=\alpha_p$, providing an alternative description
of $\frkm(P)$ in terms of meromorphic differentials on $X$ with poles ``no worse than $\frkm(P)$''; see Proposition \ref{prop:alphap-diff}.

As another application of our work, we deduce from Theorem \ref{thm:main-intro}
a description of the connected component of the identity of $\pi^{\Nori}(U)^{\ab}$ in terms of the Frobenius-power torsion subgroup
schemes of generalized Jacobians, which is a perhaps more satisfyingly direct analogue of \eqref{eq:Nori}
than the description of the full $\pi^{\Nori}(U)^{\ab}$ afforded by \eqref{eq:nori-affine}; see Corollary \ref{cor:conn part}.

In the final \S\ref{sec:singular curve}, as a related topic, we study the abelian unipotent fundamental group scheme of the singular curve $X_{\frkm}$ associated to a modulus $\frkm$ \cite[IV, \S1, No.\,4]{Serre},  and deduce a version of Das's description of the \'etale fundamental group of singular curves \cite[Theorem 3.3]{Das}; see Theorem \ref{thm:uni ab}.

\begin{rem}
\hspace{1em}
\begin{enumerate}
\item In \cite{CH}, Campbell and Hayash reformulate and reprove Geometric Class Field theory with arbitrary ramification by using the theory of {\it commutative group stacks}. Although they do not consider non-reduced finite abelian $k$-group schemes, it seems that our Theorem \ref{thm-intro:CFT} might also be deducible from their work. 

\item In \cite{Brion}, Brion has extended the theory of the Serre--Oort fundamental group to allow for arbitrary base fields (not necessarily algebraically closed).  
Relying on Brion's work, it would be interesting to generalize the present work to general base fields as well. 
\end{enumerate}
\end{rem}

\begin{ack}
We would like to thank H\'el\`ene Esnault for giving us fruitful advice about how to explain the status and position of Theorem \ref{thm-intro:CFT}. We would also like to thank Lei Zhang for informing us about the recent book by Lei Fu on algebraic curves and class field theory. We would like to express our gratitude to the referee for a careful reading of the first version of this paper and for many valuable comments. We also thank Kang Huang for giving us helpful comments. The first named author was supported by NSF DMS-2302072. 
The second named author was supported by JSPS KAKENHI Grant Numbers JP21K20334, JP24K16894.
\end{ack}

\section*{Notation and conventions}

\begin{itemize}[leftmargin=*]

\item Throughout, we work over a fixed algebraically closed field $k=\overline{k}$ of characteristic $p>0$. A {\it $k$-group} always means an abelian group scheme of finite type over $k$. A \textit{proalgebraic} $k$-group is a projective limit $\varprojlim G_i$ of $k$-groups $G_i$. A {\it profinite} $k$-group is a projective limit of finite $k$-groups.

\item For any proalgebraic $k$-group $G$, we denote by $G^0$ the connected component of the identity of $G$, and write $G^{\et}$ for the quotient $G/G^0$ and call it the {\it \'etale part} of $G$. 
As $k$ is algebraically closed, $G^{\et}$ is canonically isomorphic to the constant $k$-group scheme associated with the profinite group of $k$-valued points $G^{\et}(k)$.

\item A finite $k$-group is said to be \textit{\'etale} if $G=G^{\et}$. A finite $k$-group $G$ is said to be of \textit{multiplicative type} if it is isomorphic to $\prod_{i=1}^m\mu_{n_i}$ for some integers $n_1,\dots,n_m>0$. 
A finite $k$-group $G$ is \textit{unipotent} if it is obtained by successive extensions of $\Z/p\Z$ or $\alpha_p$. 
For any finite $k$-group $G$,  
there exists a canonical decomposition $G=G^{\uni}\times G^{\MT}$ into the product of a unipotent $k$-group $G^{\uni}$ and a multiplicative type $k$-group $G^{\MT}$ (cf.\ \cite[9.5 Theorem]{Wat}). 
We call $G^{\uni}$ (respectively $G^{\MT}$) the unipotent part (respectively the multiplicative part) of $G$. 
A finite $k$-group $G$ is said to be \textit{local-local} if it is connected and unipotent. For any finite $k$-group $G$, we put $G^{\LL}\Def (G^0)^{\uni}$ and call it the \textit{local-local} part of $G$. For any finite $k$-group $G$, there exists a functorial decomposition $G = G^0\times G^{\et}= G^{0,\rm m}\times G^{\rm ll} \times G^{\et}$, and the terminology introduced in this paragraph can be extended to any profinite $k$-group in an obvious manner. 
\end{itemize}

\section{Finite abelian torsors}\label{sec:abelian torsors}

We maintain the notation of \S\ref{sec:intro},
so $X$ is a smooth, connected and projective 
curve over $k$ and $U= X\setminus S$
is the complement in $X$ of a finite (possibly empty) set of points $S\subset X(k)$. 
By a {\it modulus} supported on $S$, we mean an effective divisor on $X$ of the form $\frkm=\sum_{x\in S}n_x x$. For any $k$-morphism $f\colon U\to G$ into a smooth $k$-group and any modulus $\frkm$ supported on $S$, we say that $f$ admits a modulus $\frkm$, or $\frkm$ is a modulus for $f$, if and only if $\frkm$ is a modulus for the induced map 
$f\colon U(k)\to G(k)$ in the sense of \cite[III, \S1, No.\,1, Definition 1]{Serre}.  

We fix a base point $x_0\in U(k)$,
and for any modulus $\frkm$ supported on $S$, 
we write $J_{X,\frkm}$ for the generalized Jacobian
of $X$ with modulus $\frkm$, and 
$\varphi_{\frkm} : U \rightarrow J_{X,\frkm}$
for the Albanese map corresponding to $x_0$ 
\cite[V, \S2, No.\,9, Theorem 1]{Serre}. 
Due to \cite[V, \S2, No.\,9, Theorem 2]{Serre},
this
satisfies the following universal property: 
for any smooth and connected $k$-group $H$
and any morphism of pointed $k$-schemes
$f:(U,x_0)\rightarrow (H,0)$ admitting $\frkm$
as a modulus,
there exists a unique homomorphism of
$k$-groups $\theta: J_{X,\frkm}\rightarrow H$
with $f = \theta \circ \varphi_{\frkm}$.

Let $G$ be a finite $k$-group.
For any modulus $\frkm$ with support contained in $S$,
pullback along $\varphi_{\frkm}$ yields
a homomorphism
\begin{equation}\label{eq:phim}
\varphi_{\frkm}^*\colon \Ext^1(J_{X,\frkm},G)\longrightarrow H^1_{\fppf}(U,G).
\end{equation}
These are compatible with change in $\frkm$,
so they induce a homomorphism of abelian groups
\begin{equation}\label{eq:CFT}
\varinjlim_{\frkm}\Ext^1(J_{X,\frkm},G)\longrightarrow H^1_{\fppf}(U,G),
\end{equation}
with $\frkm$ ranging over all moduli supported on $S$.

\begin{thm}\label{thm:CFT}
The map (\ref{eq:CFT}) is an isomorphism.
\end{thm}

We first prove:

\begin{lem}\label{lem:CFT inj}
The map
\eqref{eq:phim} is injective.
\end{lem}

\begin{proof}
Suppose that $0\to G\to P\xrightarrow{f} J_{X,\frkm}\to 0$ is an extension of $J_{X,\frkm}$ by $G$ whose pullback 
$\varphi_{\frkm}^*P\to U$ to a $G$-torsor over $U$ is trivial, and let 
$s\colon U\to\varphi_{\frkm}^*P$ be a corresponding section.
Composing $s$ with the canonical projection $\varphi_{\frkm}^*P\rightarrow P$ gives a lift $\varphi:U\rightarrow P$ of $\varphi_{\frkm}$. 
Let $P_{\rm red}$ be the reduction (cf.\ \cite[Definition 01J4]{SP}) of $P$; it is a closed subgroup scheme of $P$ by \cite[Lemma 047R]{SP},
and hence a smooth $k$-group.  
As $U$ is reduced, the map $\varphi\colon U\to P$ factors as $\varphi\colon U\to P_{\rm red}\subset P$. 
\[
\xymatrix{
P_{\rm red}\ar@{^{(}->}[r]&P\ar[d]^f\\
U\ar@{-->}[u]\ar[r]_{\varphi_{\frkm}}\ar[ur]^{\varphi}&J_{X,\frkm}
}
\]
By applying \cite[III, \S2, No.\,9, Proposition 14]{Serre} or \cite[Lemma 2.2.12]{Fu} to the composition map $\varphi_{\frkm}=f\circ\varphi\colon U\to P_{\rm red}\to J_{X,\frkm}$, we find that $\frkm$ is a modulus for $\varphi\colon U\to P_{\red}$, hence there exists a homomorphism $\theta\colon J_{X,\frkm}\to P_{\red}$ such that $\varphi=\theta\circ\varphi_{\frkm}+\varphi(x_0)$. As $f\circ\varphi(x_0)=\varphi_{\frkm}(x_0)=0$, we have  $\varphi_{\frkm}=f\circ\theta\circ\varphi_{\frkm}$, hence by the universal property of $\varphi_{\frkm}$, we can conclude that  $f\circ\theta={\rm id}_{J_{X,\frkm}}$. 
Thus, $\theta$ is a splitting of the given extension $P$, whose
class in $\Ext^1(J_{X,\frkm},G)$ is therefore trivial. 
\end{proof}

\begin{rem}\label{rem:CFT inj}
Let $f\colon H\to G$ be a finite surjective $k$-homomorphism of $k$-groups. If $G$ is reduced and connected, then we have the commutative diagram of short exact sequences
\[
\xymatrix{
1\ar[r]&\Ker(f|_{H^0_{\red}})\ar@{^{(}->}[d]\ar[r]&H^0_{\red}\ar[r]^{f|_{H^0_{\red}}}\ar@{^{(}->}[d]&G\ar[r]\ar@{=}[d]&1\\
1\ar[r]&\Ker(f)\ar[r]&H\ar[r]^f&G\ar[r]&1
}
\]
Clearly, the map $H^0_{\red}\hookrightarrow H$ is equivariant with respect to the actions of $\Ker(f|_{H^0_{\red}})$ and $\Ker(f)$. Thus, the diagram tells us that $H=\Ind_{\Ker(f|_{H^0_{\red}})}^{\Ker(f)}(H^0_{\red})$ as a $\Ker(f)$-torsor over $G$. 
\end{rem}

\begin{proof}[Proof of Theorem \ref{thm:CFT}]
Thanks to 
Lemma \ref{lem:CFT inj},
it remains to prove that \eqref{eq:CFT} is surjective.
Via the canonical decomposition 
$G=G^{0,\MT} \times G^{\rm ll} \times G^{\et}$,
the surjectivity of \eqref{eq:CFT} is reduced to the three cases
$G=G^{\star}$ for each value of $\star$ in the decomposition above. The case $\star={\et}$ is classical \cite[VI,\S2]{Serre}.

Assume now that $G=G^{0,\rm m}$. If $S=\emptyset$, then the surjectivity of
\eqref{eq:CFT} is again classical; see \S\ref{sec:intro}.
 Supposing that $S\neq\emptyset$, let $\frkm_{\red}\Def\sum_{x\in S}x$. As $\frkm_{\red}$ is reduced, $J_{X,\frkm_{\red}}$ is a semi-abelian variety over $k$. More precisely, there exists an exact sequence of smooth $k$-groups
\[
\xymatrix{
0\ar[r]& L_{\frkm_{\red}}\ar[r]& J_{X,\frkm_{\red}}\ar[r]& J_X\ar[r]&0,
}
\]
where $L_{\frkm_{\red}}$ is a smooth $k$-group, which is non-canonically isomorphic to $\G_{m,k}^{\oplus \#S-1}$; see Remark \ref{rem:L_m}; see also \cite[Proposition 3.4.9]{Fu} for a precise description of $L_{\frkm_{\red}}$. We will show that the injective homomorphism 
\[
\varphi_{\frkm_{\red}}^*\colon\Ext^1(J_{X,\frkm_{\red}},G^{0,\MT})\longrightarrow H^1_{\fppf}(U,G^{0,\MT})
\]
is bijective. As the source and target are finite groups, it suffices to prove that
\[
\#\Ext^1(J_{X,\frkm_{\red}},G^{0,\MT})=\#H^1_{\fppf}(U,G^{0,\MT}).
\]
Since $G^{0,\MT}=\prod_{i=1}^m\mu_{p^{n_i}}$ for some $n_1,\dots,n_m$, we are reduced to the case where $G^{0,\MT}=\mu_{p^n}$. 
By \cite[Proposition 3.2]{Otabe}, we have 
$H^1_{\fppf}(U,\mu_{p^n})\simeq(\Z/p^n\Z)^{\oplus f_X+\#S-1}$, where $f_X$ denotes the $p$-rank of $X$. 
Therefore, it suffices to show that $\Ext^1(J_{X,\frkm_{\red}},\mu_{p^n})\simeq(\Z/p^n\Z)^{\oplus f_X+\#S-1}$.   
By \cite[II, \S12, Lemma (12.2)]{Oort}, 
we have $\Ext^2(J_X,\mu_{p^n})=0$, and hence a short exact sequence
\[
\xymatrix{
0\ar[r]& \Ext^1(J_X,\mu_{p^n})\ar[r]& \Ext^1(J_{X,\frkm_{\red}},\mu_{p^n})\ar[r]&\Ext^1(L_{\frkm_{\red}},\mu_{p^n})\ar[r]& 0.
}
\]
As $\Ext^1(J_X,\mu_{p^n})\simeq 
J_X[p^n](k)\simeq(\Z/p^n\Z)^{\oplus f_X}$ (cf.\ \cite[\S6, p.61]{Mumford}; see also \ref{eq:Raynaud}), 
it suffices to show that $\Ext^1(L_{\frkm_{\red}},\mu_{p^n})\simeq(\Z/p^n\Z)^{\oplus \#S-1}$. Since $L_{\frkm_{\red}}\simeq\G_{m,k}^{\oplus \#S-1}$, this follows from the exactness of the sequence
\[
\xymatrix{
0\ar[r]& \End(\G_{m,k})\ar[r]^{p^n}&\End(\G_{m,k})\ar[r]&\Ext^1(\G_{m,k},\mu_{p^n})\ar[r]& 0.
}
\]
This proves the surjectivity of the map (\ref{eq:CFT}) when $G=G^{0,\MT}$.

It remains to show the surjectivity of (\ref{eq:CFT}) when $G=G^{\rm ll}$. If $S=\emptyset$, this is classical; see \S\ref{sec:intro}. Thus, we may assume that $U$ is affine. Thanks to \cite[V, \S1, Proposition 2.5]{DG},
we may embed $G$ as a subgroup scheme of a smooth connected unipotent $k$-group $W$ with quotient $W/G$ that is necessarily smooth as $G$ acts 
freely on $W$ \cite[Lemma 05B5]{SP}.  On the other hand, 
$H^1_{\fppf}(U,W)=0$ as $W$ is a successive extension of
copies of $\G_a$ (cf.\ \cite[\S1.3, Proposition 7(d)]{Serre_GP}, \cite[Exp.\,XVII, Definition 1.1]{SGA3ii}) and $H^1_{\fppf}(U,\G_a)=H^1(U,\scrO_U)=0$ 
by \cite[III, \S3, Proposition 3.7]{Milne} and the fact that $U$ is affine.
We thereby obtain
the exact sequence
\begin{equation*}
    \xymatrix{
        0 \ar[r] & {H^0(U,W)} \ar[r] & {H^0(U,W/G)} \ar[r] & {H^1_{\fppf}(U,G)} \ar[r] & 0
    }.
\end{equation*}
It follows that 
for any $G=G^{\rm ll}$-torsor $P\to U$, there exists a $k$-morphism $f\colon U\to W/G$ such that $P\simeq f^*W$ as $G$-torsors over $U$. Replacing $f$ by $f-f(x_0)$ does not change the class of $P$. We may therefore assume that $f(x_0)=0$.
 As $W/G$ is a smooth connected $k$-group, there exists a modulus $\frkm$ supported on $S$ such that $f\colon U\to W/G$ factors through $\varphi_{\frkm}\colon U\to J_{X,\frkm}$ (see  \cite[V, \S2, No.\,9, Theorem 2]{Serre}, \cite[Theorem 2.2.9]{Fu}). 
This proves that the class of the $G$-torsor $P\to U$ in $H^1_{\fppf}(U,G)$ lies in the image of $\varphi_{\frkm}^*\colon\Ext^1(J_{X,\frkm},G)\to H^1_{\fppf}(U,G)$, as desired.
\end{proof}

In particular, the proof of Theorem \ref{thm:CFT} together with \cite[VI, \S2, No.\,12, Examples.\,(1)]{Serre} shows:

\begin{cor}\label{cor:CFT MT}
Let $\frkm_{\red}\Def\sum_{x\in S}x$. 
If $G$ is of multiplicative type, we have an isomorphism
\[
\Ext^1(J_{X,\frkm_{\red}},G)\xrightarrow{~\simeq~} H^1_{\fppf}(U,G).
\]
\end{cor}

For an arbitrary finite $k$-group $G$, Theorem \ref{thm:CFT} yields an evident filtration on $H^1_{\fppf}(U,G)$:

\begin{definition}\label{def:F_m H^1}
For any modulus $\frkm$ supported on $S$, we define
\[
F_{\frkm}H^1_{\fppf}(U,G)\Def\im\left(\varphi_{\frkm}^*\colon\Ext^1(J_{X,\frkm},G)\to H^1_{\fppf}(U,G)\right).
\]
\end{definition}

\begin{rem}
The definition of $F_{\frkm}H^1_{\fppf}(U,G)$ does not depend on the choice of the base point $x_0\in U(k)$. Indeed, let $x_1\in U(k)$ be another $k$-point and $\varphi'_{\frkm}\colon U\to J_{X,\frkm}$ the Albanese map associated with $x_1$. According to \cite[V, \S2, No.\,10, Corollary 1]{Serre}, there exists a translation map $T\colon J_{X,\frkm}\to J_{X,\frkm}$ such that $\varphi'_{\frkm}=T\circ\varphi_{\frkm}$, and it follows
at once that $\varphi_{\frkm}^*$ and $(\varphi_{\frkm}')^*$
have the same image.
\end{rem}

The set of all moduli with support contained in $S$ admits an obvious partial ordering that makes this set a lattice 
(explicitly, we have $\frkm \le \frkm'$ if and only if the divisor $\frkm'-\frkm$
is effective); in particular, ${\rm Inf}(\frkm,\frkm')$ is defined for any two moduli with support contained in $S$.

\begin{prop}\label{pr:minmod}
Let $G$ be a finite $k$-group and $\frkm$, $\frkm'$ moduli supported on $S$. Then
\[
F_{\frkm}H^1_{\fppf}(U,G)\cap F_{\frkm'}H^1_{\fppf}(U,G)=F_{{\rm Inf}(\frkm,\frkm')}H^1_{\fppf}(U,G).
\]
\end{prop}

\begin{proof}
Using Lemma \ref{lem:CFT inj}, the proof of \cite[VI, \S2, No.\,12, Lemma 1]{Serre} (which is stated in the context of finite \'etale $k$-groups $G$)
works for any finite $k$-group $G$; see also \cite[Proposition 3.5.7]{Fu}.
\end{proof}

As a consequence, we obtain our version of \cite[VI, \S2, No.\,12, Proposition 11]{Serre}.

\begin{cor}\label{cor:minmod}
Let $P\to U$ be a torsor for a finite $k$-group $G$.
The subset of moduli $\frkm$ with support contained in $S$ for which
there exists an isogeny $J'\rightarrow J_{X,\frkm}$ with kernel $G$ such that $P$ is the pullback of $J'$ along $\varphi_{\frkm}:U\rightarrow J_{X,\frkm}$
is nonempty, and has a unique least (in the partial ordering specified above) element.
\end{cor}

\begin{definition}\label{defn:minmodulus}
 The {\em minimal modulus} of $P\to U$ is the unique least modulus guaranteed by Corollary \ref{cor:minmod}.
\end{definition}

\begin{rem}\label{rem:L_m}
As in \cite[V,\S3]{Serre} and \cite[Proposition 3.4.9]{Fu}, the structure
of generalized Jacobians can be made somewhat explicit.
Let $\frkm=\sum_{x\in S} n_x x$ be a modulus with support $T:={\rm Supp}(\frkm)\subseteq S$.  
When $\frkm=0$, we have $J_{X,0}=J_X$. On the other hand, 
when $\frkm\neq 0$, there exists an exact sequence of $k$-groups
\begin{equation}\label{eq:genjac-str}
    \xymatrix{
        1 \ar[r] & {\G_{\rm m}} \ar[r]^-{\Delta} & { \displaystyle\prod_{x\in T} \left(\G_{\rm m}\times V_{(n_x)}\right)} \ar[r] & {J_{X,\frkm}} \ar[r] & J_{X} \ar[r] & 0 
    }
\end{equation}
in which $\Delta$ is the diagonal embedding of 
$\G_{\rm m}$ into $\prod_{x\in T} \G_{\rm m}$ and $V_{(n)}$ is the smooth connected unipotent
$k$-group determined by 
$$
    V_{(n)}(k)=(1+uk[\![u]\!])^{\times} / (1+u^n k[\![u]\!])^{\times}
$$
for $n\ge 1$; in addition, put $V_{(0)}=0$, the trivial group scheme, for $n=0$. 
Note that $L_{\frkm}\Def\Ker(J_{X,\frkm}\twoheadrightarrow J_X)$
is a linear algebraic group.  
Moreover, as $k$ has characteristic $p$, the unipotent group $V_{(n)}$ is a product of finite-length Witt groups.  More precisely, 
let $F(u)\Def\exp(-\sum_{s\ge 0} u^{p^s}/p^s)\in\Q[\![u]\!]$;
as is well known, $F(u)$
has coefficients in $\Z_{(p)}$ so may be considered as
a formal power series over $\F_p$ by reduction modulo $p$.
For any Witt vector ${\bm a}=(a_0,a_1,a_2,\ldots)\in W(k)$,
set
$$
    E(\bm{a}\cdot u)\Def F(a_0u)F(a_1u^p)F(a_2u^{p^2})\cdots,
$$
which is an element of $1+uk[\![u]\!]$.
For each integer $i$ with $1\le i < n$ and $(i,p)=1$,
let $r_i$ be the least integer $r$ satisfying $p^r \ge n/i$.  Then
the map
\begin{equation}\label{eq:local grp}
\prod_{1\le i< n,~p\nmid i}W_{r_i}\to V_{(n)}
\quad\text{determined on $k$-points by}\quad (\bm{a}_i)_i\mapsto \prod_iE(\bm{a}_i\cdot u^i)
\end{equation}
is an isomorphism of $k$-groups \cite[V, \S3, No.\,16, Proposition 9]{Serre}.
\end{rem}

\begin{rem}\label{rem:I_m}
As in \cite[VI,\S6]{Serre} and \cite[\S2.2]{Fu}, it is illuminating to 
interpret the exact sequence \eqref{eq:genjac-str}
in terms of  {\it id\`eles}. 
For each closed point $x\in X(k)$, we denote by $\widehat{K}_x$ the completion of the function field $k(X)$ with respect to the (normalized) valuation $v_x$ associated with $x$, and we write $\widehat{\scrO}_{x}$ for the ring of integers of $\widehat{K}_x$ with maximal ideal $\frkm_x$ and set $U_x=U_x^{(0)}\Def\widehat{\scrO}_{x}^{\times}$.
For $n\ge 1$, put
$$U_x^{(n)}\Def\{g\in \widehat{K}_x^{\times}\colon v_x(1-g) \ge n\}=1+\frkm_x^n.$$
By definition, the  group of {\it  id\`eles} of $k(X)$ is the subgroup $I$ of $\prod_{x\in X(k)}\widehat{K}_x^{\times}$ consisting of families of elements $a=\{a_x\}_{x\in X(k)}$ such that $a_x\in U_x$ for almost all $x$. For any modulus $\frkm=\sum_{x\in S} n_xx$ with support $T={\rm Supp}(\frkm)\subseteq S$, we define the subgroup $I_{\frkm}$ of $I$ to be
\[
I_{\frkm}\Def\prod_{x\in T}U_x^{(n_x)}\times\prod_{x\not\in T}U_x;
\]
for $\frkm=0$, this means $I_0=\prod_{x\in X(k)}U_x$. 
For $\frkm\neq 0$, the exact sequence of $k$-points arising from \eqref{eq:genjac-str}
may be identified with the top row in the commutative diagram
\begin{equation*}
    \xymatrix{
        0 \ar[r] & k^{\times} \ar[r]\ar@{=}[d] & \displaystyle\prod_{x\in T} U_x/U_x^{(n_x)} \ar[r]\ar@{=}[d] & C_{\frkm}^0 \ar[r]\ar@{^{(}->}[d] & C^0 \ar[r]\ar@{^{(}->}[d] & 0 \\
        0 \ar[r] & k^{\times} \ar[r] & I_0/I_{\frkm} \ar[r] & I/I_{\frkm}k(X)^{\times} \ar[r] & I/I_0k(X)^{\times} \ar[r] & 0
    }
\end{equation*}
in which 
$$C_{\frkm}\Def\mathrm{Div}(X\setminus T)/\{(f)\,|\,f\equiv1~ \text{mod}~\frkm\}$$ 
is the group of $\frkm$-equivalence classes of divisors prime to $T$, and $C^0_{\frkm}$  is the subgroup of $C_{\frkm}$ 
consisting of degree-zero divisor classes (see \cite[VI, \S6]{Serre} and \cite[\S1.4]{Fu}). Note that $I/I_{\frkm}k(X)^{\times}$
is naturally identified with the group $C_{\frkm}$ \cite[Proposition 2.2.3]{Fu}, 
so that the cokernel of the rightmost two vertical maps
above is $\Z$ in each case, via the degree map. 

Furthermore, for any pointed map $f\colon (U,x_0)\to (H,0)$ into a smooth connected $k$-group $H$ which admits $\frkm$ as a modulus, if we denote by $(f,-)_{X,x}^H$ the corresponding {\it local symbol} in the sense of \cite[III,\S1, No.\,1, Proposition 1]{Serre} and \cite[\S2.2]{Fu} (see also \S\ref{sec:torsors with modulus}),  we obtain an associated homomorphism
\begin{equation*}
    \theta_f\colon I/k(X)^{\times}\to H(k)\quad\text{determined by}\quad
    \theta_f(a)=\sum_{x\in X(k)}(f,a_x)_{X,x}^H,
\end{equation*}
which satisfies $\theta_f(I_{\frkm})=0$ (see \cite[III, \S1, No.\,1, p.30]{Serre}). 
Via the identification $J_{X,\frkm}(k)=C_{\frkm}^0$, the induced homomorphism $\theta_f\colon I/I_{\frkm}k(X)^{\times}\to H(k)$ yields a homomorphism of abelian groups $\theta_f|_{L_{\frkm}(k)}\colon L_{\frkm}(k)\to H(k)$. 

On the other hand, if  
\[
\theta\colon J_{X,\frkm}\longrightarrow H
\] 
is the morphism associated with the map $f\colon (U,x_0)\to (H,0)$, whose existence follows from \cite[V,\S2, No.\,9, Theorem 2]{Serre}, according to \cite[V, \S3, No.\,18]{Serre}, we have
\begin{equation}\label{eq:theta_f}
\theta|_{L_{\frkm}(k)}(a)=-\sum_{x\in S}(f,a_x)^H_{X,x}=-\theta_f|_{L_{\frkm}(k)}(a)
\end{equation}
for any $a\in L_{\frkm}(k)$. 
\end{rem}

\begin{rem}\label{rem:repartition}
Let $R$ be the ring of {\it r\'epartitions} for the curve $X$ in the sense of \cite[II, No.\,5]{Serre}, i.e.\ $R$ is the subring of $\prod_{x\in X(k)}k(X)$ consisting of families $r=(r_x)_{x\in X(k)}$ of elements of $k(X)$ such that $r_x\in\scrO_{X,x}$ for almost all $x\in X(k)$, and let $R(0)\Def\prod_{x\in X(k)}\scrO_{X,x}$ be the $k$-subalgebra of $R$ consisting of
those families $r$ with $r_x \in \scrO_{X,x}$ for all $x$. 
By \cite[II, No.\,5, Proposition 3]{Serre}, we have a natural identification
\[
H^1_{\fppf}(X,\G_a)=R/(R(0)+k(X)).
\]
One can extend this to an analogous identification
with an arbitrary smooth connected unipotent $k$-group $H$ in place of $\G_a$ as follows. 
We define $H(R)$ to be the subgroup of $\prod_{x\in X(k)}H(k(X))$ consisting of families $h=(h_x)_{x\in X(k)}$ of elements of $H(k(X))$ such that $h_x\in H(\scrO_{X,x})$ for almost all $x\in X(k)$. We also consider the subgroup $H(R(0))\Def\prod_{x\in X(k)}H(\scrO_{X,x})\subset H(R)$. 
We claim that there is a natural identification of abelian groups
\begin{equation}\label{eq:repartition-H1}
H^1_{\fppf}(X,H)\simeq H(R)/(H(R(0))+H(k(X))).
\end{equation}
To see this, let $H_X$ be the Zariski sheaf on $X$ determined by the assignment $H_X(W)\Def H(W)$ for any open subset $W\subseteq X$. 
As $H$ is a successive extension of copies of $\G_a$, we have $H^1_{\fppf}(K,H)=0$ for any field $K$ over $k$, hence the natural identification
\[
H^1_{\fppf}(X,H)\simeq H^1(X,H_X).
\]
Consider the short exact sequence of Zariski sheaves on $X$:
\[
\xymatrix{
0\ar[r]&H_X\ar[r]&H(k(X))\ar[r]&H(k(X))/H_X\ar[r]&0,
}
\]
where $H(k(X))$ is considered as a constant sheaf on $X$, hence a flasque sheaf, which guarantees the vanishing of the first cohomology $H^1(X,H(k(X)))$. 
This yields an exact sequence
\[
\xymatrix{
H(k(X))\ar[r]&H^0(X,H(k(X))/H_X)\ar[r]&H^1(X,H_X)\ar[r]&0.
}
\]
For any $s\in H(k(X))$, there exists an open subset $W\subseteq X$ such that $s$ belongs to the image of the restriction map $H_X(W)=H(W)\to H(k(X))$, and we have natural identifications
\[
H^0(X,H(k(X))/H_X)\simeq\bigoplus_{x\in X(k)}H(k(X))/H(\scrO_{X,x})\simeq H(R)/H(R(0))
\]
as in the proof of \cite[II, No.\,5, Proposition 3]{Serre}. This establishes the identification (\ref{eq:repartition-H1}). 
\end{rem}

\section{Finite local-local torsors over affine curves}\label{sec:torsors with modulus}

As an illustration of the previous considerations,
we now work out explicitly the classification of finite local-local torsors in terms of Corollary \ref{cor:minmod}, and its relation
to other well-known classifications. 
Let $X$ be a smooth, connected and projective curve over $k$ and $H$ a smooth connected $k$-group. For any rational map $f\in H(k(X))$, we consider the associated {\it local symbol}
\begin{equation}\label{eq:loc symbol for X}
\left\{(f,-)^H_{X,x}\colon k(X)^{\times}\longrightarrow H(k)\right\}_{x\in X(k)}
\end{equation}
in the sense of \cite[III, \S1]{Serre}. 
More precisely, if we take a finite subset $S\subset X(k)$  so that $f$ defines a $k$-morphism $f\colon X\setminus S\to H$, then
by Rosenlicht’s theorem  \cite[III, \S1, No.\,3, Theorem 1]{Serre} there exists a modulus $\frkm$ for $f$ supported on $S$,
and it follows from \cite[III, \S1, No.\,1, Proposition 1]{Serre} that there is a unique {\em local symbol} in the sense
of \cite[III, \S1, No.\,1, Definition 2]{Serre} associated to $f$ and $\frkm$.  By \cite[III, \S1, No.\,1, Remark]{Serre}, this local symbol depends only
on $f$. 
Throughout this section, we use the notation introduced in Remarks \ref{rem:L_m} and \ref{rem:I_m}.

Fix $x\in X(k)$, and consider the local symbol at $x$: 
\[
(-,-)^H_{X,x}\colon H(k(X))\times k(X)^{\times}\longrightarrow H(k).
\]
The local symbol $(-,-)_{X,x}^H$ at $x$ factors through the composition of natural maps: 
\begin{align*}
H(k(X))\times k(X)^{\times}&\xrightarrow{~(a)~} H(k(X))/H(\scrO_{X,x},\frkm_x)\times k(X)^{\times}\\
&\xrightarrow{~(b)~} H(k(X))/H(\scrO_{X,x},\frkm_x)\times \varprojlim_nk(X)^{\times}/U_{X,x}^{(n)}\\
&\xrightarrow[~\simeq~]{~(c)~}H(\widehat{K}_x)/H(\widehat{\scrO}_{x},\frkm_x)\times\varprojlim_nk(X)^{\times}/U_{X,x}^{(n)}\\
&\xrightarrow[~\simeq~]{~(d)~} H(\widehat{K}_x)/H(\widehat{\scrO}_{x},\frkm_x)\times\varprojlim_{n}\widehat{K}_x^{\times}/U_x^{(n)},
\end{align*}
where we set $U_{X,x}^{(0)}\Def \scrO_{X,x}^{\times}$ and  $U_{X,x}^{(n)}\Def \{g\in\scrO_{X,x}^{\times}\,|\,v_x(1-g)\ge n\}$ for $n\ge 1$; 
\begin{align*}
H(\scrO_{X,x},\frkm_x)&\Def\Ker\left(H(\scrO_{X,x})\to H(k)\right),\quad H(\widehat{\scrO}_{x},\frkm_x)\Def\Ker\left(H(\widehat{\scrO}_x)\to H(k)\right).
\end{align*}
The first map $(a)$ is induced by properties (ii), (iii) of \cite[III, \S1, No.\,1, Definition 2]{Serre}. Indeed, these properties imply that for any $f\in H(\scrO_{X,x},\frkm_x)\subseteq H(\scrO_{X,x})$, we have $(f,-)^H_{X,x}=v_x(-)f(x)=0$. Moreover, for a fixed $f\in H(k(X))$, 
a modulus condition, together with property (i) of \cite[III, \S1, No.\,1, Definition 2]{Serre}, supplies an integer $n\ge 0$ such that the map $(f,-)^H_{X,x}\colon k(X)^{\times}\to H(k)$ factors through the natural homomorphism $k(X)^{\times}\twoheadrightarrow k(X)^{\times}/U_{X,x}^{(n)}$, hence we obtain the map $(b)$; the last identification $(d)$ follows from the natural isomorphism
\[
k(X)^{\times}/U_{X,x}^{(n)}\xrightarrow{~\simeq~}\widehat{K}_x^{\times}/U_x^{(n)}
\]
for any $n\ge 0$; 
the bijectivity of the third map $(c)$ follows from the natural identification
\[
H(k(X))/H(\scrO_{X,x},\frkm_x)
\xrightarrow{~\simeq~}H(\widehat{K}_x)/H(\widehat{\scrO}_{x},\frkm_x).
\]
To see this identification, recall that $H$ fits into an
exact sequence
\[
0\longrightarrow B\longrightarrow H
\longrightarrow A\longrightarrow 0,
\]
where $A$ is an abelian variety and $B$ is a smooth connected affine $k$-group. 
By the properness of $A$, we have $A(k(X))=A(\scrO_{X,x})$, $A(\widehat{K}_x)=A(\widehat{\scrO}_x)$, hence 
\[
A(k(X))/A(\scrO_{X,x},\frkm_x)=A(k)=A(\widehat{K}_x)/A(\widehat{\scrO}_x,\frkm_x).
\] 
Moreover, since $k$ is algebraically closed, $B$ is a
successive extension of copies of $\G_a$ and $\G_m$.
In particular, the vanishing of the first cohomology groups
with coefficients in $\G_a$ and $\G_m$ over fields and
local rings implies that $H^1_{\fppf}(S,B)=0$ for $S=k(X),\scrO_{X,x},\widehat{K}_x,\widehat{\scrO}_x$. 
Consequently, we have the following commutative diagram
with exact rows:
\[
\xymatrix{
0\ar[r]&
B(k(X))/B(\scrO_{X,x},\frkm_x)
\ar[r]\ar[d]&
H(k(X))/H(\scrO_{X,x},\frkm_x)
\ar[r]\ar[d]&
A(k)\ar[r]\ar@{=}[d]&0\\
0\ar[r]&
B(\widehat{K}_x)/B(\widehat{\scrO}_x,\frkm_x)
\ar[r]&
H(\widehat{K}_x)/H(\widehat{\scrO}_x,\frkm_x)
\ar[r]&
A(k)\ar[r]&0
}
\]
Thus, we are reduced to the case $H=B$. In this case, $H$ is a product of a split torus and a
smooth connected unipotent $k$-group.
For $H=\G_m$, we have the natural identifications
\[
\G_m(k(X))/\G_m(\scrO_{X,x},\frkm_x)
\simeq k^\times\times\Z
\simeq
\G_m(\widehat{K}_x)/
\G_m(\widehat{\scrO}_x,\frkm_x).
\]
If $H$ is a smooth connected unipotent $k$-group,
then, since $k$ is algebraically closed, it is a
successive extension of copies of $\G_a$.
Using the vanishing of the first cohomology groups
with coefficients in $\G_a$ on affine schemes,
we are reduced to the case $H=\G_a$.
In this case, every element of $\widehat{K}_x$
is congruent modulo $\frkm_x$ to a rational function,
obtained by retaining the finite principal part and
the constant term of its Laurent expansion. Here, the Laurent expansion
must be in terms of a uniformizing parameter in $\scrO_{X,x}$, rather than $\widehat{\scrO}_{x}$. 
Hence the natural map
\[
\G_a(k(X))/\G_a(\scrO_{X,x},\frkm_x)
\xrightarrow{~\simeq~}
\G_a(\widehat{K}_x)/
\G_a(\widehat{\scrO}_x,\frkm_x)
\]
is an isomorphism.

We denote by
\[
\overline{(-,-)}_{X,x}^H\colon H(\widehat{K}_x)
\times\widehat{K}_x^{\times}\longrightarrow H(\widehat{K}_x)/H(\widehat{\scrO}_{x},\frkm_x)\times\varprojlim_{n}\widehat{K}_x^{\times}/U_x^{(n)}\longrightarrow H(k)
\]
the induced bilinear map. Namely, 
given $\widehat f\in H(\widehat K_x)$, choose $f\in H(k(X))$
representing its class modulo $H(\widehat{\scrO}_x,\frkm_x)$. Choose $N\ge1$ such that $(f,U_{X,x}^{(N)})^H_{X,x}=0$.
For $\widehat g\in\widehat K_x^\times$, choose $g\in k(X)^\times$
with $\widehat g/g\in U_x^{(N)}$, and set
\[
\overline{(\widehat f,\widehat g)}^H_{X,x}=(f,g)^H_{X,x}.
\]
Changing $g$ does not change this value by the modulus condition.
Changing $f$ changes it by an element of
$H(\scrO_{X,x},\frkm_x)$, whose local symbol is zero.
Taking a common sufficiently large $N$ also proves independence
of $N$ and bilinearity.

\begin{lem}\label{lem:loc symbol L}
Let $L\Def k(\!(u)\!)$. 
For any $k$-isomorphism $\sigma\colon L\xrightarrow{\simeq}\widehat{K}_x$, the induced bilinear map
\[
(-,-)^H_{L}\colon H(L)\times L^{\times}\overset{\sigma}{\simeq} H(\widehat{K}_x)
\times\widehat{K}_x^{\times}\xrightarrow{~\overline{(-,-)}^H_{X,x}~} H(k)
\]
does not depend on the choice of $\sigma$, $X$, or $x\in X(k)$. 
\end{lem}

\begin{proof}
This fact appears to be well known to experts, although
we have not been able to find a precise reference in
the literature.
Let $W\subseteq X$ be an open neighborhood of $x$,
and let $w\colon W\to\bbP^1_k$ be an \'etale morphism
with $w(x)=\infty$.
For any $f\in H(k(\bbP^1_k))$ and
$g\in k(\bbP^1_k)^\times$, we have
\[
(f,g)_{\bbP^1_k,\infty}^H
=
(w^*f,w^*g)_{X,x}^H.
\]
This follows from the proof of \cite[Lemma 4.38]{RS}.

By the comparison between rational and completed
arguments established above, this equality extends
to arbitrary elements of the completed local fields. We regard $u$ as a rational uniformizer at $\infty$
on $\bbP^1_k$, so that $\widehat{\scrO}_{\bbP^1_k,\infty}=k[[u]]$. Let
\[
u\longmapsto t\Def a_1u+a_2u^2+\cdots\in k[[u]],\quad a_1\ne 0,
\]
be a formal change of local parameter. For any $m\geq 1$, its truncation
\[
t_m\Def a_1u+\cdots+a_mu^m\in k[u]
\]
satisfies
\[
t\equiv t_m\pmod{u^{m+1}}.
\]
Since $t_m'(0)=a_1\ne 0$, the rational map
$w_m\colon\bbP^1_k\dashrightarrow\bbP^1_k$
defined by $w_m^*(u)=t_m$ is an \'etale
morphism in a neighborhood of $\infty$. 
For fixed $f\in H(k(\bbP^1_k))$ and
$g\in k(\bbP^1_k)^\times$, choosing $m$
sufficiently large gives
\[
f(t)\equiv f(t_m)
\mod~ H(k[[u]],uk[[u]])
\]
and
\[
g(t)/g(t_m)\in U^{(n)}_{\infty},
\]
where $n$ is chosen according to the modulus
condition for the first variable $f$.
Consequently, the \'etale compatibility for
rational arguments extends to arbitrary formal
changes of parameter.
This proves the independence of $\sigma$, $X$,
and $x$.
\end{proof}

In the remainder of this section, we denote by $L\Def k(\!(u)\!)$ the field of Laurent formal power series over $k$ 
with valuation ring $\scrO_L\Def k[\![u]\!]$ and normalized valuation $v_L$. 
We set $U_L\Def U_L^{(0)}\Def\scrO_L^{\times}$, 
and for $n\ge 1$, define $U_L^{(n)}\Def\{g\in\scrO_L^{\times}\,|\,v_L(1-g)\ge n\}$. 
For any smooth connected
$k$-group $H$, we denote by 
\[
(-,-)^H_L\colon H(L)\times L^{\times}\longrightarrow H(k)
\] 
the {\it local symbol} for $H$ in the sense of Lemma \ref{lem:loc symbol L}.

\begin{definition}\label{def:RoSe}
Let $H$ be a smooth connected $k$-group. 
We define an increasing filtration $\{\fil_nH(L)\}_{n\ge 0}$ on $H(L)$ by 
\[
\fil_nH(L)\Def
\begin{cases}
H(\scrO_L)& \text{if $n=0$},\\
\{f\in H(L)\,|\,(f,U_L^{(n)})^H_L=0\}& \text{if $n>0$},
\end{cases}
\] 
which we call the {\it Rosenlicht--Serre filtration} on $H(L)$.
\end{definition}

\begin{rem}\label{rem:RoSe=KR}
For a finite-length Witt group $H=W_m$, an explicit description of the local symbol 
\[
(-,-)_L^{W_m}\colon W_m(L)\times L^{\times}\longrightarrow W_m(k)
\] 
is established in \cite[Proof of Proposition 6.4]{KR}.
If $\{\fil_nW_m(L)\}_{n\ge 0}$ is the Rosenlicht--Serre filtration of Definition \ref{def:RoSe} and  $\{\fil_{n-1}^FW_m(L)\}_{n\ge 1}$
is the Kato--Russell filtration of \cite[\S2]{KR}, then we have $\fil_nW_m(L) = \fil_{n-1}^FW_m(L)$ for $n\ge 1$
while $\fil_0W_m(L)=W_m(\scrO_L)=\fil_0^FW_m(L)$; see \cite[Proposition 6.4(3)]{KR}, \cite[\S7]{RS}. 
\end{rem}

\begin{prop}\label{prop:RoSe}
With $L$ as above, let $n$ be a nonnegative integer.
\begin{enumerate}
\item \label{prop:RoSe-1} For any homomorphism $\phi\colon H_1\to H_2$ of smooth connected $k$-groups, we have $\phi(\fil_nH_1(L))\subseteq\fil_nH_2(L)$.

\item \label{prop:RoSe-2} For any homomorphism $\phi\colon H_1\to H_2$ of smooth connected $k$-groups whose kernel is finite and connected, we have $\phi(\fil_nH_1(L))=\phi(H_1(L))\cap\fil_nH_2(L)$.
\end{enumerate}
\end{prop}

\begin{proof}
By \cite[III, \S1, No.\,2, Proposition 2]{Serre}, we have $(\phi(f),-)_L^{H_2}=\phi((f,-)_L^{H_1})$ for $f\in H_1(L)$,
which gives \eqref{prop:RoSe-1} for $n>0$. For $n=0$, claim \eqref{prop:RoSe-1} says that $\phi(H_1(\scrO_L))\subseteq H_2(\scrO_L)$. This is evident. Let us show \eqref{prop:RoSe-2}. By \eqref{prop:RoSe-1}, it follows that we have an inclusion $\phi(\fil_nH_1(L))\subseteq \phi(H_1(L))\cap\fil_nH_2(L)$. The reverse inclusion for $n>0$ is a straightforward consequence of the injectivity of the map 
$\phi\colon H_1(k)\to H_2(k)$
induced by $\phi$ on $k$-valued points. For $n=0$, as the map $\phi\colon H_1\to H_2$ is proper, the valuative criterion for properness implies that $\phi(H_1(L))\cap H_2(\scrO_L)\subseteq\phi(H_1(\scrO_L))$. This completes the proof. 
\end{proof}

\begin{prop}\label{prop:local symbol}
Let $H$ be a smooth connected unipotent $k$-group. 
There exists an isomorphism 
\[
H(L)/H(\scrO_L)\xrightarrow{~\simeq~}\varinjlim_n\Hom(V_{(n)},H)
\]
which maps each $\fil_nH(L)/\fil_0H(L)$ onto $\Hom(V_{(n)},H)$, and the induced homomorphism
\[
H(L)/H(\scrO_L)\longrightarrow\varinjlim_n\Hom(V_{(n)}(k),H(k))
\] 
coincides with the map given by $-(-,-)_{L}^H$, the additive  inverse of the local symbol $(-,-)_L^H\colon H(L)\times L^{\times}\to H(k)$. 
\end{prop}

\begin{proof}
Let $\A^1_k=\Spec k[t]$ with $t=1/u$ and $X\Def\bbP^1_k$ the smooth compactification of $\A^1_k$. Note that the restriction map
\[
H(\A^1_k)\longrightarrow H(L)
\]
induces an isomorphism
\[
\Hom((\A^1_k,0),(H,0))\xrightarrow{~\simeq~}H(L)/H(\scrO_L),
\]
where the left hand side means the group of pointed $k$-morphisms $(\A^1_k,0)\to (H,0)$. 
Indeed, as $H^1_{\fppf}(Z,H)=0$ for any smooth connected unipotent $k$-group $H$ and any affine $k$-scheme $Z$, d\'evissage for smooth connected unipotent $k$-groups reduces the claimed identification to the case $H=\G_a$, which is obvious. It follows that there is a unique isomorphism 
\[
H(L)/H(\scrO_L)\xrightarrow{~\simeq~}\varinjlim_n\Hom(V_{(n)},H)
\]
making the diagram
\[ \xymatrix{ \Hom((\A_k^1,0),(H,0))\ar[r]^{\xi~}_{\simeq~}\ar[d]_{\simeq}&\varinjlim_{n}\Hom(J_{X,n\cdot\infty},H)\ar@{=}[d]\\
H(L)/H(\scrO_L)\ar[r]^{\simeq\quad}&\varinjlim_n\Hom(V_{(n)},H),
} 
\] 
commute, where the upper horizontal isomorphism $\xi$ is due to \cite[V,\S2, No.\,9, Theorem 2]{Serre}. For the other statements, see the equation (\ref{eq:theta_f}) in Remark \ref{rem:I_m}. 
\end{proof}

\begin{lem}\label{lem:inj F_* Ext}
For integers $l,m,r\ge 1$, the $r$-th iterated Frobenius homomorphism $F^r\colon W_m\to W_m$ induces an injective homomorphism
\[
(F^r)_*\colon \Ext^1(W_l,W_m)\hookrightarrow\Ext^1(W_l,W_m).
\]
Consequently, the connecting homomorphism $\Hom(W_{l},W_m)\to\Ext^1(W_l,W_m[F^r])$
is surjective, where $W_m[F^r]\Def\Ker(F^r\colon W_m\to W_m)$.
\end{lem}

\begin{proof}
As $(F^r)_*=(F_*)^r$, we are reduced to the case of $r=1$. Let $W\Def\varprojlim_nW_n$ be the additive group of Witt vectors.  
By \cite[II, \S9, Corollary 9.4]{Oort}, we have
$\Ext^1(W,\G_a)=0$. 
It then follows by induction on $m$, using the exact sequence $0\to W_{m-1}\xrightarrow{V}W_m
\to\G_a\to 0$, that $\Ext^1(W,W_m)=0$ for all $m\geq 1$.
Therefore, the exact sequence
\[
\xymatrix{
0\ar[r]& W\ar[r]^{V_W^l}&W\ar[r]& W_l\ar[r]& 0,
}
\]
where $V_W^l\colon W\hookrightarrow W$ is the $l$-th iterated Verschiebung operator on Witt vectors, yields the long exact sequence
\[
\xymatrix{
0\ar[r]&\Hom(W_l,W_m)\ar[r]&\Hom(W,W_m)\ar[r]^{~(V_W^l)^*~}&\Hom(W,W_m)\ar[r]&\Ext^1(W_l,W_m)\ar[r]& 0. 
}
\]
Hence, we have 
\begin{equation}\label{eq:inj F_* Ext}
\Ext^1(W_l,W_m)\simeq\Hom(W,W_m)/(V_W^l)^*(\Hom(W,W_m)).
\end{equation}

Next, we claim that, as subgroups of $\Hom(W,W_m)$,
\begin{equation}
	(V_W^l)^*(\Hom(W,W_m)) = \begin{cases}
		(V^l)_*(\Hom(W,W_{m-l})) & \text{if}\ l < m \\
		0 & \text{if}\ l \ge m
	\end{cases},
	\label{eq:homcases}
\end{equation}
in which $(V^l)_*$ is the map induced by post-composition with the $l$-th iterated Verschiebung operator $V^l\colon W_{m-l}\hookrightarrow W_m$ on finite Witt vectors. 
Let $v_m:W_m\to W_m$ be the endomorphism
\[
v_m(a_0,\ldots,a_{m-1})=(0,a_0,\ldots,a_{m-2}).
\]
For $h:W\to W_m$, let $h^{(p^l)}$ be obtained by raising
its coordinate-polynomial coefficients to their $p^l$-th powers.
The identity $F^l h=h^{(p^l)}F_W^l$ and the Witt relations give
\[
hV_W^lF_W^l=[p^l]_{W_m}h
=v_m^lF^lh=v_m^lh^{(p^l)}F_W^l.
\]
Cancellation of the faithfully flat source Frobenius $F_W^l$ yields
$hV_W^l=v_m^lh^{(p^l)}$.
For $l\ge m$, we have $v_m^l=0$, proving the vanishing
in \eqref{eq:homcases}.
For $l<m$, write $v_m^l=V^lR^l$, where
$R^l:W_m\to W_{m-l}$ is truncation. Its kernel is isomorphic
to $W_l$, so $\Ext^1(W,W_l)=0$ implies that
\[
(R^l)_*:\Hom(W,W_m)\twoheadrightarrow\Hom(W,W_{m-l}).
\]
As $k$ is perfect, the coefficient twist is bijective.
Consequently,
\[
(V_W^l)^*\Hom(W,W_m)
=(V^l)_*\Hom(W,W_{m-l}),
\]
which proves \eqref{eq:homcases}.

Now when $l < m$, the short exact sequence
\[
\xymatrix{
0\ar[r]& W_{m-l}\ar[r]^{V^l}&W_m\ar[r]& W_l\ar[r]& 0,
}
\]
and the vanishing $\Ext^1(W,W_{m-l})=0$ yield a natural identification
\begin{equation}
\Hom(W,W_m)/(V^l)_*(\Hom(W,W_{m-l}))\simeq\Hom(W,W_l).
\label{eq:llessm}
\end{equation}
Combining \eqref{eq:llessm}, \eqref{eq:homcases}, and \eqref{eq:inj F_* Ext},
we deduce a natural identification
\begin{equation}
	\Ext^1(W_l,W_m) \simeq \Hom(W,W_{\min\{m,l\}}),
\end{equation}
which is induced by target truncation, and the map $F_*\colon\Ext^1(W_l,W_m)\to\Ext^1(W_l,W_m)$ then becomes post-composition by $F\colon W_{\min\{m,l\}}\to W_{\min\{m,l\}}$ on $\Hom(W,W_{\min\{m,l\}})$. 
Hence, the claimed injectivity follows from the vanishing 
\[
\Hom(W,W_{\min\{m,l\}}[F])=0.
\]
This is valid because $W$ is reduced, whereas $W_{\min\{m,l\}}[F]$ is finite
and infinitesimal, so every homomorphism from $W$ to $W_{\min\{m,l\}}[F]$ is zero.
\end{proof}

\begin{rem}\label{rem:inj F_* Ext}
Let $G=G^{\LL}$ be a finite local-local $k$-group.  Let $0\to G\to H_1\to H_2\to 0$ be a resolution of $G$ by smooth connected unipotent $k$-groups $H_j$,  
and let $p^n$ be the {\it period} of $H_1$, i.e.\ $p^n$ is the smallest power $p^m$ of $p$ satisfying $p^mH_1=0$ \cite[VII, \S2, No.\,10]{Serre}. 
Note that $H_2$ must then also have period $p^n$. 
For $l\ge n$, we claim that the induced map
\[
\Ext^1(W_l,H_1)\longrightarrow\Ext^1(W_l,H_2)
\]
is injective. 

Indeed, since $p^l=0$ on $H_j$, 
the map
\[
(V_W^l)^*\colon \Hom(W,H_j)\longrightarrow\Hom(W,H_j)
\]
induced by the Verschiebung operator $V_W^l$ on $W$ is the zero map. 
To see this, recall that $p^l=V_W^l\circ F_W^l$ on $W$.
Hence $0=(p^l)^*=(F_W^l)^*\circ(V_W^l)^*$ on $\Hom(W,H_j)$. Since $F_W^l\colon W\to W$ is faithfully flat,
the induced map $(F_W^l)^*$ is injective, and therefore $(V_W^l)^*=0$. 
Since each $H_j$ has a filtration with successive quotients $\G_a$, the vanishing $\Ext^1(W,\G_a)=0$ implies 
$\Ext^1(W,H_j)=0$ by induction. 
This together with the exact sequence
\[
0\longrightarrow W\xrightarrow{~V^l~}W
\longrightarrow W_l\longrightarrow 0
\]
yields a natural
isomorphism
\[
\Hom(W,H_j)\xrightarrow{~\simeq~}\Ext^1(W_l,H_j).
\]
The claimed injectivity then follows from
$\Hom(W,G)=0$. 
\end{rem}

\begin{definition}\label{def:adm-resol}
Let $G=G^{\LL}$ be a finite local-local $k$-group. A resolution $0\to G\to H_1\to H_2\to 0$ of $G$ by smooth connected unipotent $k$-groups $H_j$ is said to be \textit{admissible} if for any finite-length Witt group $W_l$, the connecting homomorphism $\Hom(W_l,H_2)\to\Ext^1(W_l,G)$ is surjective. 
\end{definition}

\begin{rem}\label{rem:adm-resol}
\hspace{1em}
\begin{enumerate}
\item By Lemma \ref{lem:inj F_* Ext}, for any $m,r>0$, the exact sequence
\[
\xymatrix{
0\ar[r]& W_m[F^r]\ar[r]& W_m\ar[r]^{F^r}&W_m\ar[r]& 0
}
\]
is an admissible resolution of $W_m[F^r]$. 

\item Any finite local-local $k$-group $G=G^{\LL}$ admits an admissible resolution. 
Indeed, for sufficiently large $m, N$, there exists an embedding $G\hookrightarrow W_m^N$, and as $G$ is local, we must have 
$G\subset W_m[F^r]^N$ for some $r>0$. This inclusion extends to a commutative diagram with exact rows
\[
\xymatrix{
0\ar[r]&G\ar[r]\ar@{^{(}->}[d]&W_m^N\ar[r]\ar@{=}[d]&W_m^N/G\ar[r]\ar@{->>}[d]&0\\
0\ar[r]&W_m[F^r]^N\ar[r]&W_m^N\ar[r]^{F^r}&W_m^N\ar[r]&0
}
\]
For any $l>0$, this yields the commutative diagram
\[
\xymatrix{ 
\Ext^1(W_l,W_m^N)\ar@{=}[d]\ar[r]&\Ext^1(W_l,W_m^N/G)\ar[d]&\\ 
\Ext^1(W_l,W_m^N)\ar@{^{(}->}[r]^{F^r_*}&\Ext^1(W_l,W_m^N)
} 
\]
and it follows that the resolution $0\to G\to W_m^N\to W_m^N/G\to 0$ of $G$ is admissible. 

\item Given a commutative diagram of resolutions of $G=G^{\LL}$ by smooth connected unipotent $k$-groups
\[
\xymatrix{
0\ar[r]&G\ar[r]\ar@{=}[d]&H_1\ar[r]\ar[d]&H_2\ar[r]\ar[d]&0\\
0\ar[r]&G\ar[r]&H'_1\ar[r]&H'_2\ar[r]&0
}
\]
If the upper resolution is admissible, then so is the lower one. This is immediate from the commutativity of the associated square
\[
\xymatrix{ \Hom(W_l,H_2)\ar@{->>}[r]\ar[d]&\Ext^1(W_l,G)\ar@{=}[d]\\ \Hom(W_l,H'_2)\ar[r]&\Ext^1(W_l,G)
} 
\]
\end{enumerate}
\end{rem}

\begin{prop}\label{prop:conductor vs adm resol}
Recall that $L\Def k(\!(u)\!)$ and let
$$ ({\rm R})\qquad 0\longrightarrow G\longrightarrow H_1\longrightarrow H_2\longrightarrow 0$$
be an admissible resolution of $G=G^{\LL}$. For any $n\ge 0$, 
we denote by $\fil_n^{({\rm R})}H^1_{\fppf}(L,G)$ the image of $\fil_nH_2(L)$ under the associated connecting homomorphism $H_2(L)\twoheadrightarrow H^1_{\fppf}(L,G)$. Then there exists a natural isomorphism
\begin{equation}\label{eq:prop:conductor vs adm resol}\dfrac{\fil_n^{({\rm R})}H^1_{\fppf}(L,G)}{\fil_0^{({\rm R})}H^1_{\fppf}(L,G)}\xrightarrow{~\simeq~}\Ext^1(V_{(n)},G)
\end{equation}
for any $n\ge 0$, where both sides are zero for $n=0$, 
and 
the increasing filtration $\{\fil_n^{({\rm R})}H^1_{\fppf}(L,G)\}_{n\ge 0}$ on $H^1_{\fppf}(L,G)$ is independent of the choice of $({\rm R})$. \end{prop}

\begin{proof}
First, note that the filtration at level zero is independent of the choice of admissible resolution. More precisely, for every admissible resolution $({\rm R})$, we have
\[
\fil_0^{({\rm R})}H^1_{\fppf}(L,G)=\im(H^1_{\fppf}(\scrO_L,G)\longrightarrow H^1_{\fppf}(L,G)).
\]
Indeed, the connecting map over $\scrO_L$ is surjective because $H^1_{\fppf}(\scrO_L,H_1)=0$, and it commutes with restriction to $L$.

Fix $n> 0$. 
As $H^{1}_{\fppf}(Z,H_j)=0$ for any affine $k$-scheme $Z$, we have the exact sequence
\begin{equation}\label{eq:H2toH1fppfsurj}
0\longrightarrow\dfrac{H_1(L)}{H_1(\scrO_L)}\longrightarrow\dfrac{H_2(L)}{H_2(\scrO_L)} \longrightarrow\dfrac{H^1_{\fppf}(L,G)}{H^1_{\fppf}(\scrO_L,G)} \longrightarrow 0.
\end{equation}
By Proposition \ref{prop:RoSe}, 
we thus obtain the exact sequence
\[
0\longrightarrow\dfrac{\fil_nH_1(L)}{\fil_0H_1(L)}\longrightarrow\dfrac{\fil_nH_2(L)}{\fil_0H_2(L)}\longrightarrow\dfrac{\fil_n^{({\rm R})}H^1_{\fppf}(L,G)}{\fil_0^{({\rm R})}H^1_{\fppf}(L,G)}\longrightarrow 0.
\]
By Proposition \ref{prop:local symbol}, the additive inverses of the local symbols $\fil_nH_j(L)\to\Hom(V_{(n)},H_j)$ for $j=1,2$ make the following diagram with exact rows commute:
\[
\xymatrix{
0\ar[r]&\dfrac{\fil_nH_1(L)}{\fil_0H_1(L)}\ar[r]\ar[d]^{\simeq}&\dfrac{\fil_nH_2(L)}{\fil_0H_2(L)}\ar[r]\ar[d]^{\simeq}& \dfrac{\fil_n^{({\rm R})}H^1_{\fppf}(L,G)}{\fil_0^{({\rm R})}H^1_{\fppf}(L,G)}\ar[r]& 0\\
0\ar[r]&\Hom(V_{(n)},H_1)\ar[r]&\Hom(V_{(n)},H_2)\ar[r]&\Ext^1(V_{(n)},G)\ar[r]&0
}
\]
which yields the claimed isomorphism (\ref{eq:prop:conductor vs adm resol}). 

For the independence of the choice of the admissible resolution $({\rm R})$, let
\[
({\rm R}')\qquad
0\longrightarrow G \longrightarrow H'_1
\longrightarrow H'_2\longrightarrow 0
\]
be another admissible resolution of $G$. 
Let $i\colon G\hookrightarrow H_1$ and $i'\colon G\hookrightarrow H'_1$ denote the corresponding
embeddings. Consider the diagonal embedding
\[
(i,i')\colon G\hookrightarrow H_1\times H'_1.
\]
Since both $({\rm R})$ and $({\rm R}')$ are admissible,
the induced maps
\[
i_*\colon \Ext^1(W_l,G)\longrightarrow
\Ext^1(W_l,H_1),\quad i'_*\colon  \Ext^1(W_l,G)\longrightarrow\Ext^1(W_l,H'_1)
\]
are zero for any $l\ge 1$. By the additivity of $\Ext^1(W_l,-)$, the map
\[
(i,i')_*:
\Ext^1(W_l,G)\longrightarrow
\Ext^1(W_l,H_1\times H'_1)
\]
is also zero for any $l\ge 1$. Thus, the long exact sequence associated with
\[
0\longrightarrow G\longrightarrow H_1\times H'_1
\longrightarrow (H_1\times H'_1)/G\longrightarrow 0
\]
shows that this resolution is admissible. 
We therefore reduce to the case where there exists a commutative diagram of admissible resolutions:
\[
\xymatrix{
0\ar[r]&G\ar[r]\ar@{=}[d]&H'_1\ar[r]\ar[d]&H'_2\ar[r]\ar[d]&0\\
0\ar[r]&G\ar[r]&H_1\ar[r]&H_2\ar[r]&0.
}
\]
In this case, the isomorphism \eqref{eq:prop:conductor vs adm resol} implies that the induced map
\[
\dfrac{\fil_n^{({\rm R'})}H^1_{\fppf}(L,G)}{\fil_0^{({\rm R'})}H^1_{\fppf}(L,G)}\longrightarrow \dfrac{\fil_n^{({\rm R})}H^1_{\fppf}(L,G)}{\fil_0^{({\rm R})}H^1_{\fppf}(L,G)}
\]
is bijective. Note that the identification \eqref{eq:prop:conductor vs adm resol} is natural for morphisms of admissible resolutions and the induced map on $\Ext^1(V_{(n)},G)$ is the identity. 
Since
$\fil_0^{({\rm R'})}H^1_{\fppf}(L,G)=\fil_0^{({\rm R})}H^1_{\fppf}(L,G)$ and $\fil_n^{({\rm R'})}H^1_{\fppf}(L,G)
\subseteq\fil_n^{({\rm R})}H^1_{\fppf}(L,G)$,
it therefore follows that
\[
\fil_n^{({\rm R'})}H^1_{\fppf}(L,G)
=
\fil_n^{({\rm R})}H^1_{\fppf}(L,G).
\]
This establishes the claimed independence of the choice of
admissible resolution and completes the proof.
\end{proof}

\begin{definition}\label{def:cond}
Let $G=G^{\LL}$ be a finite local-local $k$-group. 
We define an increasing filtration $\{\fil_nH^1_{\fppf}(L,G)\}_{n\ge 0}$ by
\[
\fil_nH^1_{\fppf}(L,G)\Def\im(\fil_nH_2(L)\subset H_2(L)\twoheadrightarrow H^1_{\fppf}(L,G))
\]
for any admissible resolution $0\to G\to H_1\to H_2\to 0$ of $G$; by Proposition \ref{prop:conductor vs adm resol}, this is independent of the choice of admissible resolution.
For any $n\ge 0$, we denote by $\psi^L_n$ the isomorphism
\[
\dfrac{\fil_nH^1_{\fppf}(L,G)}{\fil_0H^1_{\fppf}(L,G)}\xrightarrow{~\simeq~}\Ext^1(V_{(n)},G)
\]
of Proposition \ref{prop:conductor vs adm resol} and put $\varphi^L_n\Def(\psi^L_n)^{-1}$. 
We denote by $\psi^L$ the limit of $\psi^L_n$, i.e. 
\[  
\psi^L\colon\dfrac{H^1_{\fppf}(L,G)}{H_{\fppf}^1(\scrO_L,G)}=\varinjlim_n\dfrac{\fil_nH^1_{\fppf}(L,G)}{\fil_0H^1_{\fppf}(L,G)}\xrightarrow{~\simeq~}\varinjlim_n\Ext^1(V_{(n)},G).
\]
Furthermore, we define a map $\cond_L^{G}\colon H^1_{\fppf}(L,G)\to\Z_{\ge 0}$ by
\[
\cond_L^{G}(P)\Def\min\{n\ge 0\,|\,P\in\fil_nH^1_{\fppf}(L,G)\}.
\]
\end{definition}

\begin{rem}\label{rem:RS}
In \cite[\S9]{RS}, R\"ulling and Saito use the theory of
{\em reciprocity sheaves} to provide a new conductor for torsors
for arbitrary finite $k$-groups $G$. 
Assume that $K$ is a geometric henselian discrete valuation field over $k$, and let
$$
    c_{K}^{H^1(G)}: H^1_{\fppf}(K,G) \longrightarrow \Z_{\ge 0}
$$
be the {\em motivic conductor} of \cite{RS} associated to the presheaf with transfers $H^1(G)\Def H^1_{\fppf}(-,G)$. 
On the other hand, when $G=\alpha_p$ and $L=k(\!(u)\!)$, using the admissible resolution $0\to\alpha_p\to\G_a\xrightarrow{F}\G_a\to 0$ 
of $\alpha_p$, Proposition \ref{prop:conductor vs adm resol} implies
\[
\cond_L^{\alpha_p}(P)=\min\{n\ge 0\,|\,P\in\im(\fil_n\G_a(L)\to H^1_{\fppf}(L,\alpha_p))\}
\]
for any $P\in H^1_{\fppf}(L,\alpha_p)$. 

To compare these two conductors, let $K=\Frac(\scrO^h_{X,x})$ and $L=\widehat{K}_x\simeq k((u))$.
Since $k$ is perfect, the explicit filtration of
\cite[\S9.9, Proposition 9.10]{RS} and the Rosenlicht--Serre
filtration for $H=\G_a$ have the same images in $H^1_{\fppf}(L,\alpha_p)=L/L^p$: the image of
$\scrO_L$ at level zero and of $u^{1-n}\scrO_L$ at level $n\ge1$.
Thus the motivic conductor extends to $L$ by
\[
\widehat c_L^{H^1(\alpha_p)}([a])=
\begin{cases}
0,&da\in\Omega^1_{\scrO_L/k},\\
-\operatorname{ord}_x(da),&\text{otherwise}.
\end{cases}
\]
Indeed, after removing negative Laurent terms that are $p$-th
powers, a largest remaining pole of order $j$ gives conductor
$j+1$. This description is independent of the parameter.
The same calculation over $K$, together with preservation of
differential orders under completion, gives
\[
c_K^{H^1(\alpha_p)}(P)
=\widehat c_L^{H^1(\alpha_p)}(P|_L),
\qquad
\widehat c_L^{H^1(\alpha_p)}=\cond_L^{\alpha_p}.
\]
\end{rem}

\begin{thm}\label{thm:m vs c}
Let $G=G^{\LL}$ be a finite local-local $k$-group and $P\in H^1_{\fppf}(U,G)$ have minimal modulus $\frkm(P)=\sum_{x\in S}n_x(P) x$ in the sense of Definition \ref{defn:minmodulus}.
 For any $x\in S$, we have
\[
n_x(P)=\cond^G_{\widehat{K}_{x}}(P|_{\widehat{K}_{x}}).
\]
\end{thm}

\begin{lem}\label{lem:m vs c}
For any finite local-local $k$-group $G$, the map $\varphi_0^*\colon\Ext^j(J_X,G)\to H^j_{\fppf}(X,G)$ is bijective for $j=1,2$. 
\end{lem}

\begin{proof}
For $j=1$, this is well known. We will show the claim for $j=2$.  
As $J_X$ is an abelian variety, we have $\Ext^3(J_X,H)=0$ for any $k$-group $H$ by \cite[II,\S12, Proposition (12.5)]{Oort}, so 
a d\'evissage reduces us to the case $G=\alpha_p$. 
Thanks to the vanishing  $H^2_{\fppf}(X,\G_a)=H^2_{\rm Zar}(X,\scrO_X)=0$ and $\Ext^2(J_X,\G_a)=0$ \cite[II,\S12, Lemma (12.8)]{Oort}, we have a commutative diagram with exact rows 
\begin{equation*}
    \xymatrix{
    0\ar[r]&\Ext^1(J_X,\alpha_p)\ar[r]\ar[d]&\Ext^1(J_X,\G_a)\ar[r]\ar[d]&\Ext^1(J_X,\G_a)\ar[r]\ar[d]&\Ext^2(J_X,\alpha_p)\ar[r]\ar[d]& 0\\
        0\ar[r]& H^1_{\fppf}(X,\alpha_p)\ar[r]& H^1_{\fppf}(X,\G_a)\ar[r]& H^1_{\fppf}(X,\G_a)\ar[r]& H^2_{\fppf}(X,\alpha_p)\ar[r]& 0
    }
\end{equation*}
whose vertical maps are induced by pullback along an Albanese map $\varphi_0\colon X\rightarrow J_X$.
The first vertical map is an isomorphism by the $j=1$ instance of Lemma \ref{lem:m vs c},
while the second and third vertical maps are the composites of the isomorphisms 
$\Ext^1(J_X,\G_a)\simeq H^1(J_X,\scrO_{J_X})\simeq H^1(X,\scrO_X)$ (see \cite[VII, \S3--\S4, Theorems 7, 8 and 9]{Serre}) with the identification $H^i(X,\scrO_X)=H^i_{\fppf}(X,\G_a)$ for $i=1$. 
The final vertical map must therefore be an isomorphism as well.
\end{proof}

\begin{lem}\label{lem:m vs c 2}
Suppose that $S\neq\emptyset$. For any smooth connected unipotent $k$-group $H$, we have 
$\varinjlim_{\frkm}\Ext^1(J_{X,\frkm},H)=0$.
\end{lem}

\begin{proof}
Let $\frkm$ be a modulus supported on $S$ and $0\to H\to J\to J_{X,\frkm}\to 0$ an extension of smooth $k$-groups. As $S\neq\emptyset$, we have $H^1_{\fppf}(U,H)=0$ by the affineness of $U=X\setminus S$, and the Albanese map $\varphi_{\frkm}\colon U\to J_{X,\frkm}$ lifts to a map $f:U\rightarrow J$. Replacing $f$ with $f-f(x_0)$, by \cite[III, \S1, No.3, Theorem 1]{Serre} and \cite[V,§2, No.9, Theorem 2]{Serre}, 
$f$ admits a modulus $\frkm' \ge \frkm$ supported on $S$,
and there is a unique homomorphism $\theta: J_{X,\frkm'}\rightarrow J$
factoring $f$.
Then $\theta$ is a section of the pullback of
$J$ along $J_{X,\frkm'}\twoheadrightarrow J_{X,\frkm}$,
so the class of $J$ in $\Ext^1(J_{X,\frkm'},H)$ is trivial. This proves the lemma. 
\end{proof}

\begin{lem}\label{lem:m vs c 3}
Let $H$ be a smooth connected unipotent $k$-group. Under the identification (\ref{eq:repartition-H1}):
\[
H^1_{\fppf}(X,H)\simeq H(R)/(H(R(0))+H(k(X)))
\] 
given in Remark \ref{rem:repartition}, for any $x\in X(k)$, the following diagram is commutative
\begin{equation}\label{eq:repartition}
\xymatrix{
H(k(X))/H(\scrO_{X,x})\ar[d]_{(-,-)^H_{X,x}}\ar[r]^{\quad\delta_x}&H^1_{\fppf}(X,H)\\
\varinjlim_n\Hom(L_{n\cdot x},H)\ar[r]^{\quad\eta_x}&\Ext^1(J_X,H),\ar[u]_{\varphi_0^*}
}
\end{equation}
where the upper horizontal arrow $\delta_x$ is the composition of the natural maps
\[
H(k(X))/H(\scrO_{X,x})\hookrightarrow\bigoplus_{y\in X(k)}H(k(X))/H(\scrO_{X,y})\simeq H(R)/H(R(0))\twoheadrightarrow H^1_{\fppf}(X,H),
\]  
and the lower one $\eta_x$ sends each homomorphism $\alpha\colon L_{n\cdot x}\to H$ to the pushout $\alpha_*J_{X,n\cdot x}$ of the extension $0\to L_{n\cdot x}\to J_{X,n\cdot x}\to J_X\to 0$. 
\end{lem}

\begin{proof}
To see this, let $f\in H(k(X))\setminus H(\scrO_{X,x})$ and suppose that $\frkm$ is a modulus for $f$. As $f\not\in H(\scrO_{X,x})$, we have $x\in{\rm Supp}(\frkm)$. We put $\frkm_0:=n_x\cdot x$, where $n_x:=v_x(\frkm)>0$ is the order of $\frkm$ at $x$. Choosing $x_0\in U(k)$ and replacing $f$ with $f-f(x_0)$ (which does not change the
class of $f$ in $H(k(X))/H(\scrO_{X,x})$), the universal property of the Albanese map $\varphi_{\frkm}\colon U:=X\setminus{\rm Supp}(\frkm)\to J_{X,\frkm}$
associated to $x_0$ implies that there exists a homomorphism of $k$-groups $\theta\colon J_{X,\frkm}\to H$ such that $f=\theta\circ\varphi_{\frkm}$. We write $\theta_{x}$ for the restriction of $\theta$ to the $k$-subgroup $L_{\frkm_0}=V_{(n_x)}$ of $J_{X,\frkm}$, where the canonical map $J_{X,\frkm}\twoheadrightarrow J_{X,\frkm_0}$ maps the subgroup $V_{(n_x)}$ of $J_{X,\frkm}$ isomorphically onto the subgroup $L_{\frkm_0}$ of $J_{X,\frkm_0}$, and we identify $L_{\frkm_0}$ with the subgroup $V_{(n_x)}$ of $J_{X,\frkm}$. Then we have 
\[
\theta_{x}=-(f,-)_{X,x}^H\in\Hom(L_{\frkm_0},H)
\]
(see (\ref{eq:theta_f}) in Remark \ref{rem:I_m}). Therefore, it suffices to show that 
\begin{equation}\label{eq:lem:m vs c 3}
-\varphi^*_0(\eta_x(\theta_{x}))=\delta_x(f).
\end{equation}

According to \cite[VII,\S4, No.\,20]{Serre}, there exists $\psi_x\in L_{\frkm_0}(k(X))$ with $\varphi_{\frkm_0}-\psi_x\in J_{X,\frkm_0}(\scrO_{X,x})$ and $-(\psi_x,-)^{L_{\frkm_0}}_{X,x}={\rm id}_{L_{\frkm_0}}$. By functoriality of local symbols,
$$-(\theta_{x}\circ\psi_x,-)^H_{X,x}= \theta_{x}\circ\bigl(-(\psi_x,-)^{L_{\mathfrak m_0}}_{X,x}\bigr) = \theta_{x} = -(f,-)^H_{X,x}.$$ 
It then follows from Proposition \ref{prop:local symbol} that $f=\theta_{x}\circ\psi_x$ in $H(k(X))/H(\mathcal O_{X,x}).$ 
This implies that $\delta_x(f)=(\theta_{x})_*\delta_x(\psi_x)$. On the other hand, we have $\varphi^*_0(\eta_x(\theta_{x}))=\varphi_0^*((\theta_{x})_*J_{X,\frkm_0})=(\theta_{x})_*\varphi_0^*(J_{X,\frkm_0})$. As a consequence, to prove the identification (\ref{eq:lem:m vs c 3}), it suffices to show that  
\begin{equation}\label{eq:lem:m vs c 3-2}
\varphi^*_0(J_{X,\frkm_0})=-\delta_x(\psi_x). 
\end{equation}
We have the commutative diagram: 
\[
\xymatrix{
0\ar[r]&L_{\frkm_0}\ar[r]&J_{X,\frkm_0}\ar[r]&J_X\ar[r]&0\\
&L_{\frkm_0}\ar@{=}[u]\ar@{^{(}->}[r]&\varphi_0^*J_{X,\frkm_0}\ar[r]\ar[u]&X\ar[u]^{\varphi_0}&
}
\]

Since $L_{\frkm_0}$ is a smooth connected unipotent $k$-group, the $L_{\frkm_0}$-torsor $J_{X,\frkm_0}\rightarrow J_X$ is Zariski locally trivial. 
Let $\{U_i\}$ be a Zariski open covering of $J_X$ with $\varphi_0(x)\in U_{i_0}$ such that the morphism $J_{X,\frkm_0}\to J_X$ admits a section $s_i$ over $U_i$  for each $i$. The $f_{ij}\Def s_j-s_i\in L_{\frkm_0}(U_i\cap U_j)$ form a $1$-cocycle with values in $L_{\frkm_0}$, whose class in $\Ext^1(J_X,L_{\frkm_0})$ is represented by $J_{X,\frkm_0}$. On the other hand, the $g_{ij}:=f_{ij}\circ\varphi_0\in L_{\frkm_0}(\varphi_0^{-1}(U_i)\cap\varphi_0^{-1}(U_j))$ form a $1$-cocycle on $X$ whose class is represented by $\varphi_0^*(J_{X,\frkm_0})$. 
Note that the composition $X\setminus\{x\}\hookrightarrow X\xrightarrow{\varphi_0}J_X$ factors through the morphism $\varphi_{\frkm_0}\colon X\setminus\{x\}\to J_{X,\frkm_0}$ by the universal property of $\varphi_{\frkm_0}$. 
For any $i$, we define $h_i\in L_{\frkm_0}(k(X))$ to be
\[
h_i\Def \varphi_{\frkm_0}-s_i\circ\varphi_0.
\] 
As in the proof of \cite[VII, \S4, No.\,19, Proposition 15]{Serre}, the $h_i$ define an element $h=(h_y)_{y\in X(k)}\in L_{\frkm_0}(R)/L_{\frkm_0}(R(0))$ such that its image in $H^1_{\fppf}(X,L_{\frkm_0})$ is equal to the class of $-\varphi_0^*(J_{X,\frkm_0})$. More precisely, for $y\in X(k)$, choose $i$ with $\varphi_0(y)\in U_i$, and let $h_y$ be the class of $h_i$ in
	$L_{\frkm_0}(k(X))/L_{\frkm_0}(\scrO_{X,y}).$ Since $h_j-h_i=-(s_j-s_i)\circ\varphi_0=-g_{ij}$
	is regular wherever both terms are defined, $h_y$ is independent of $i$. 
Moreover, as 
\[
h_{i_0}-\psi_x=(\varphi_{\frkm_0}-s_{i_0}\circ\varphi_0)-\psi_x=(\varphi_{\frkm_0}-\psi_x)-s_{i_0}\circ\varphi_0\in L_{\frkm_0}(\scrO_{X,x}),
\] 
we have $h_x=\psi_x$ in $L_{\frkm_0}(k(X))/L_{\frkm_0}(\scrO_{X,x})$. On the other hand, as $\varphi_{\frkm_0}$ is regular on $X\setminus\{x\}$, we have $h_y=0$ in $L_{\frkm_0}(k(X))/L_{\frkm_0}(\scrO_{X,y})$ for $y\neq x$. This yields the desired identification (\ref{eq:lem:m vs c 3-2}).  
\end{proof}

\begin{proof}[Proof of Theorem \ref{thm:m vs c}]
If $S=\emptyset$, the assertion is vacuous. We therefore assume
$S\ne\emptyset$ throughout the proof. 
Let $H^r_S(X,G)$ (respectively $H^r_c(U,G)$) be the fppf cohomology with supports (respectively the fppf cohomology with compact supports) in the sense of \cite[III, \S0]{Milne_ADT}, \cite{DH}. We have a commutative diagram with exact rows
\begin{equation}\label{eq:proof of thm:m vs c 0}
\xymatrix{
0\ar[r]&H^1_c(U,G)\ar[r]
\ar@{^{(}->}[d]&H^1_{\fppf}(U,G)\ar[r]\ar@{=}[d]&\prod_{x\in S}H^1_{\fppf}(\widehat{K}_{x},G)\ar[r]\ar@{->>}[d]&H^2_c(U,G)
\ar@{->>}[d]\\
0\ar[r]&H^1_{\fppf}(X,G)\ar[r]\ar@{=}[d]&H^1_{\fppf}(U,G)\ar[r]\ar@{=}[d]&H^2_S(X,G)\ar@{->>}[r]\ar[d]_-{\simeq}&H^2_{\fppf}(X,G)\ar@{=}[d]\\
0\ar[r]&H^1_{\fppf}(X,G)\ar[r]&H^1_{\fppf}(U,G)\ar[r]&\prod_{x\in S}\dfrac{H^1_{\fppf}(\widehat{K}_{x},G)}{H^1_{\fppf}(\widehat{\scrO}_{x},G)}
&H^2_{\fppf}(X,G)
}
\end{equation}
where the injectivity $H^1_c(U,G)\hookrightarrow H^1_{\fppf}(U,G)$ is thanks to the vanishing $H^0_{\fppf}(\widehat{K}_{x},G)=0$ and the surjectivity $H^2_S(X,G)\twoheadrightarrow H^2_{\fppf}(X,G)$ is due to the vanishing $H^2_{\fppf}(U,G)=0$. 
Let 
\[
\xymatrix{
0\ar[r]& G\ar[r]& H_1\ar[r]& H_2\ar[r]& 0
}
\] 
be an admissible resolution of $G$ by smooth connected unipotent $k$-groups $H_j$. 
Thanks to Lemma \ref{lem:m vs c 2}, the connecting homomorphisms for variable $\frkm$ induce a surjection
\begin{equation}\label{eq:proof of thm:m vs c: surj1}
\xymatrix{
\varinjlim_{\frkm}\Hom(J_{X,\frkm},H_2)\ar@{->>}[r]&\varinjlim_{\frkm}\Ext^1(J_{X,\frkm},G),
}
\end{equation}
while the exact sequence \eqref{eq:H2toH1fppfsurj} gives a surjection
\begin{equation}\label{eq:proof of thm:m vs c: surj2}
\xymatrix{
H_2(\widehat{K}_{x})/H_2(\widehat{\scrO}_{x})\ar@{->>}[r]& H^1_{\fppf}(\widehat{K}_{x},G)/H^1_{\fppf}(\widehat{\scrO}_{x},G).
}
\end{equation}
On the other hand, the diagram
\begin{equation}\label{eq:proof of thm:m vs c}
\xymatrix{
H_2(U)\ar[r]\ar@{}[dr]|{\displaystyle\bigstar}&\prod_{x\in S}H_2(\widehat{K}_{x})/H_2(\widehat{\scrO}_{x})\ar[d]_{\simeq}^{(-(-,-)^{H_2}_{X,x})_{x\in S}}\ar[r]\ar@{}[dr]|{\qquad\displaystyle\blacklozenge}&H^1_{\fppf}(X,H_2)\\
\varinjlim_{\frkm}\Hom(J_{X,\frkm},H_2)\ar[u]^{\varinjlim\varphi_{\frkm}^*}\ar@{^{(}->}[r]&\prod_{x\in S}\varinjlim_{n}\Hom(V_{(n)},H_2)\ar[r]&\Ext^1(J_X,H_2)\ar[u]_{-\varphi_0^*}
}
\end{equation}
is commutative. 
Here the commutativity of the left square ($\bigstar$) (respectively the right square $(\displaystyle\blacklozenge)$) is due to equation (\ref{eq:theta_f}) in Remark \ref{rem:I_m} (respectively due to the commutativity of the diagram (\ref{eq:repartition}) in Lemma \ref{lem:m vs c 3}). 
The commutative diagrams (\ref{eq:proof of thm:m vs c 0}) and  (\ref{eq:proof of thm:m vs c}) make the following diagram with exact rows commute:
\[
\xymatrix{
0\ar[r]&H^1_{\fppf}(X,G)\ar[r]\ar@{=}[d]&H^1_{\fppf}(U,G)\ar[r]\ar@{=}[d]&H^2_S(X,G)\ar[r]\ar[d]_-{\simeq}&H^2_{\fppf}(X,G)\ar@{=}[d]\\
0\ar[r]&H^1_{\fppf}(X,G)\ar[r]&H^1_{\fppf}(U,G)\ar[r]\ar[r]\ar@{}[dr]|{\displaystyle\bigstar\bigstar}&\prod_{x\in S}\dfrac{H^1_{\fppf}(\widehat{K}_{x},G)}{H^1_{\fppf}(\widehat{\scrO}_{x},G)}\ar[r]\ar[d]_{\simeq}^{(\psi^{\widehat{K}_{x}})_{x\in S}}\ar@{}[dr]|{\qquad\displaystyle\blacklozenge\blacklozenge}&H^2_{\fppf}(X,G)\\
0\ar[r]&\Ext^1(J_X,G)\ar[u]_{\varphi_0^*}^{\simeq}\ar[r]&\varinjlim\Ext^1(J_{X,\frkm},G)\ar[u]_{\varinjlim\varphi_{\frkm}^*}^{\simeq}\ar[r]&\prod_{x\in S}\varinjlim\Ext^1(V_{(n)},G)\ar[r]&\Ext^2(J_X,G)\ar[u]_{-\varphi_0^*}^{\simeq}
}
\]
where $\psi^{\widehat{K}_x}$ is defined in Definition \ref{def:cond}. 
Here, the commutativity of the square ($\displaystyle\bigstar\bigstar$) (respectively the square ($\displaystyle\blacklozenge\blacklozenge$)) is reduced via the surjective map (\ref{eq:proof of thm:m vs c: surj1}) (respectively (\ref{eq:proof of thm:m vs c: surj2})) to the commutativity of ($\displaystyle\bigstar$) (respectively ($\displaystyle\blacklozenge$)) in the diagram (\ref{eq:proof of thm:m vs c}).  
As a consequence, for any modulus $\frkm=\sum_{x\in S}n_x x$ supported on $S$, we have a commutative diagram with exact rows
\[
\xymatrix{
0\ar[r]&H^1_{\fppf}(X,G)\ar[r]&H^1_{\fppf}(U,G)\ar@{}[dr]|{\displaystyle\square}\ar[r]&\prod_{x\in S}\dfrac{H^1_{\fppf}(\widehat{K}_{x},G)}{H^1_{\fppf}(\widehat{\scrO}_{x},G)}
\ar[r]&H^2_{\fppf}(X,G)\\
0\ar[r]&\Ext^1(J_X,G)\ar[r]\ar[u]^{\simeq}_{\varphi_0^*}&\Ext^1(J_{X,\frkm},G)\ar[r]\ar@{^{(}->}[u]_{\varphi_{\frkm}^*}&\prod_{x\in S}\Ext^1(V_{(n_x)},G)\ar@{^{(}->}[u]_{(\varphi^{\widehat{K}_{x}}_{n_x})_{x\in S}}\ar[r]&\Ext^2(J_X,G)\ar[u]^{\simeq}_{-\varphi_0^*}
}
\]
in which the first and last vertical maps are isomorphisms due to 
Lemma \ref{lem:m vs c}, where $\varphi_{n_x}^{\widehat{K}_x}$ is defined in Definition \ref{def:cond}. 
It follows that the middle square is Cartesian. 
Thus $P\in F_{\frkm}H^1_{\fppf}(U,G)$ if and only if
$\cond^G_{\widehat K_x}(P|_{\widehat K_x})\le n_x$
for every $x\in S$. Taking the least such modulus proves
the asserted equality of its coefficients with the local conductors. 
This completes the proof. 
\end{proof}

We conclude the present section by giving an explicit description of the minimal modulus $\frkm(P)$ for $P\in H^1_{\fppf}(U,\alpha_p)$ in terms of differential forms. 

Let $Z$ be any smooth $k$-variety and write $\Omega^1_{Z/k,\cl}:=\Ker(d:\Omega^1_{Z/k}\rightarrow \Omega^2_{Z/k})$ for the sheaf of 
closed differential 1-forms on $Z$ over $k$.
The {\em Cartier Operator} is an additive and $p^{-1}$-linear
map of sheaves $V: \Omega^1_{Z/k,\cl}\rightarrow \Omega^1_{Z/k}$ defined as in \cite{Seshadri}, and we write $H^0(Z,\Omega^1_{Z/k,\cl})[V]$ for the kernel of the induced map on global sections.


\begin{lem}\label{lem:fppf-diff}
    For any smooth $k$-variety $Z$, there is a functorial identification of abelian groups
    $$
        H^1_{\fppf}(Z,\alpha_p) \simeq H^0(Z,\Omega^1_{Z/k,\cl})[V].
    $$
  \end{lem}

\begin{proof}
    See \cite[III, \S4, Proposition 4.14]{Milne}.
\end{proof}

Now let $H$ be a smooth $k$-group, and denote by $m: H\times_k H\rightarrow H$ the multiplication map, and
by $\pr_i : H\times_k H \rightarrow H$ projection onto the $i$-th factor.  
If $G$ is any fppf sheaf of abelian groups, the set of {\em primitive elements} in $H^1_{\fppf}(H,G)$
is by definition the kernel of the map
\begin{equation*}
    m^* - \pr_1^*-\pr_2^* : H^1_{\fppf}(H,G)\longrightarrow H^1_{\fppf}(H\times_k H, G).
\end{equation*}

\begin{prop}\label{prop:ext-primitive}
    Let $H$ be a smooth $k$-group and $G$ a finite local $k$-group.
    The canonical map
    \begin{equation}
        \xymatrix{
            \Ext^1(H,G) \ar[r] & H^1_{\fppf}(H,G)
            }\label{eq:ext-fppf}
    \end{equation}
    is injective, with image equal to the set of primitive elements.
\end{prop}

\begin{proof}
    We first note that \eqref{eq:ext-fppf} is injective: if $E$ is an extension of $H$ by $G$
    that admits a section, then $E\simeq H\times G$ as $k$-schemes.  As $G$ is finite local
    and $H$ is smooth, the induced homomorphism $E_{\red}\rightarrow E\rightarrow H$ must be an 
    isomorphism, and it follows that
    $E \simeq H\times G$ as $k$-groups.
    Thanks to \cite[Lemma 2.4]{SOS}, the image of \eqref{eq:ext-fppf} is contained
    in the set of primitive elements.  To see that the image is exactly the set of
    primitive elements, let $P$ be a $G$-torsor of $H$ which represents a primitive element in $H^1_{\fppf}(H,G)$. We need to show that $P$ has a $k$-group structure such that the morphism $P\to H$ is a $k$-homomorphism whose kernel is isomorphic to $G$. As $m^*P=\pr_1^*P+\pr_2^*P$ in $H^1_{\fppf}(H\times_k H,G)$, we have a commutative diagram
    \[
    \xymatrix{
    P\times_kP\ar[r]^{\quad m_P}\ar[d]&P\ar[d]\\
    H\times_k H\ar[r]^{\quad m}&H,
    }
    \]
    where the upper horizontal arrow, denoted by $m_P\colon P\times_kP\to P$, is compatible with the actions of $G\times_k G$ and $G$ with respect to the multiplication $G\times G\to G$ on $G$ as in the proof of \cite[Theorem 2.6]{SOS}. Moreover, as $P\to H$ is finite and radicial, the induced map on $k$-points $P(k)\to H(k)$ is bijective. Hence, there is a unique $e_P \in P(k)$ lifting the identity element $e \in H(k)$,
which must moreover satisfy $m_P(e_P\times e_P)=e_P$. 
The proof of \cite[VII, \S3, No.\,15, Theorem 5]{Serre} cited there 
carries over to our situation {\em mutatis mutandis} due to the fact that any map from a smooth group such as $H$ or $H\times_k H$ to the local $G$ must be constant, from which it follows that the maps $m_P$ and $e_P$ endow $P$ with a $k$-group structure. It is also clear from our construction that the morphism $P\to H$ is a morphism of $k$-groups whose kernel coincides with $G$. This completes the proof. 
\end{proof}


\begin{cor}\label{cor:ext-diff}
    Let $H$ be a smooth $k$-group, and denote by $\omega_H$
    the subspace of $H^0(H,\Omega^1_{H/k})$ consisting of translation invariant 1-forms.
    There is a functorial isomorphism of abelian groups
    $$
        \Ext^1(H,\alpha_p) \simeq \omega_H[V]
    $$
\end{cor}

\begin{proof}
    The functorial isomorphism of Lemma \ref{lem:fppf-diff} and 
    Proposition \ref{prop:ext-primitive} yield a natural identification
    of $\Ext^1(H,\alpha_p)$ with the kernel of $m^*-\pr_1^*-\pr_2^*$
    on $H^0(H,\Omega^1_{H/k,\cl})[V]$.  Since every translation invariant 1-form on $H$ is closed \cite[Lemma 2.1]{Coleman},
    it therefore suffices to observe that $\omega_H$ is precisely the kernel of 
    $m^*-\pr_1^*-\pr_2^*$ on $H^0(H,\Omega^1_{H/k})$, which follows from \cite[III,\S3,No.\,11, Proposition 17]{Serre} applied to $(f,g)=(\pr_1,\pr_2)$. 
\end{proof}

\begin{prop}\label{prop:alphap-diff}
    Let $P\in H^1_{\fppf}(U,\alpha_p)$ and let $\omega(P)\in H^0(U,\Omega^1_X)[V]$ be the corresponding differential form.
        Then the minimal modulus $\frkm(P)$ of $P$ is given by
        $$
        \frkm(P)=\sum_{x\in S} \max\{0,-{\rm ord}_x(\omega(P))\} \cdot x
        $$
\end{prop}

\begin{proof}
Pullback along $\varphi_{\frkm}:U\rightarrow J_{X,\frkm}$ induces 
a commutative diagram
\begin{equation*}
\xymatrix{
{\omega_{J_{X,\frkm}}[V]}\ar[dd]^-{\simeq}_-{\varphi_{\frkm}^*} \ar[rr]^-{\simeq} & & \Ext^1(J_{X,\frkm},\alpha_p) \ar[d]_-{\varphi_{\frkm}^*}^-{\simeq}\\
  & & F_{\frkm}H^1_{\fppf}(U,\alpha_p) \ar@{^{(}->}[d]\\
H^0(X,\Omega^1_{X}(\frkm))[V] \ar@{^{(}->}[r] & H^0(U,\Omega^1_X)[V] \ar[r]^-{\simeq}& H^1_{\fppf}(U,\alpha_p),
}
\end{equation*}
in which $H^0(X,\Omega^1_{X}(\frkm))$ denotes the vector space of meromorphic differential forms $\eta$ on $X$ such that $\mathrm{div}(\eta)+\frkm\ge 0$, 
and the horizontal isomorphisms are induced by the identification of Lemma \ref{lem:fppf-diff},
via Corollary \ref{cor:ext-diff}, and the left vertical isomorphism is due to \cite[V,\S2, No.\,10, Proposition 5]{Serre}.
It follows that the identification $H^0(U,\Omega^1_{X/k})[V]\simeq H^1_{\fppf}(U,\alpha_p)$
carries the subspace of differential forms with poles ``no worse'' than $\frkm$
isomorphically onto $F_{\frkm}H^1_{\fppf}(U,\alpha_p)$, as claimed.
\end{proof}

\section{The abelian fundamental group scheme}\label{sec:ab Nori}

We will continue to use the same notation as in the previous sections. For a fixed base point $x_0\in U(k)$, we denote by $\pi^{\Nori}(U,x_0)$ the {\it fundamental group scheme} for $U$ in the sense of \cite[Chapter II]{Nori82}. This is a profinite $k$-group which classifies pointed $G$-torsors $(P,p)\to (U,x_0)$ for any finite (not necessarily abelian) $k$-group scheme $G$. If $x_1\in U(k)$ is any other choice of base point, then $\pi^{\Nori}(U,x_1)$ is an inner twist of $\pi^{\Nori}(U,x_0)$ \cite[II, Proposition 4]{Nori82};
in particular, the maximal abelian quotient of $\pi^{\Nori}(U,x_0)$ does not depend on the choice of the base point $x_0$, and we denote it simply by $\pi^{\ab}(U)$. This is a profinite $k$-group with the property that there exists a natural isomorphism of abelian groups 
\[
\Hom(\pi^{\ab}(U),G)\xrightarrow{~\simeq~}H^1_{\fppf}(U,G)
\]
for variable finite $k$-groups $G$.

As an application of Theorem \ref{thm:CFT}, we will provide a description of the profinite $k$-group $\pi^{\ab}(U)$ in terms of {\it Serre--Oort fundamental groups} \cite[II, \S7]{Oort}, \cite{Brion} of generalized Jacobians. Let $\calC$ be the category of $k$-groups (not necessarily connected or reduced); it is an {\it artinian} abelian category \cite[II, \S6, Lemma (6.1)]{Oort}.  
Let $\Pro(\calC)$ be the associated {\it pro-category} of $\calC$ in the sense of \cite[I, \S4]{Oort}, i.e.\ the objects of $\Pro(\calC)$ are projective systems $\{G_i\}$ of $k$-groups $G_i\in\calC$, and for any two objects $\{G_i\},\{H_j\}\in\Pro(\calC)$, the group of morphisms is given by
\[
\Hom(\{G_i\},\{H_j\})\Def\varprojlim_j\varinjlim_i\Hom(G_i,H_j).
\]
Let $\Hom(\calC,\Ab)$ be the abelian category of left exact covariant functors from $\calC$ to the category of abelian groups $\Ab$. It turns out that the functor
\[
\Pro(\calC)\longrightarrow\Hom(\calC,\Ab)^{\rm op};~\{G_i\}\longmapsto \varinjlim_i\Hom(G_i,-)
\]
is an equivalence of categories, and  
the category $\Pro(\calC)$ is, in turn, an abelian category with enough projective objects \cite[I, \S4]{Oort}. 
Let $\calF\subset\calC$ be the full subcategory of $\calC$ consisting of finite $k$-groups and $\Pro(\calF)$ the associated pro-category. 
The functor $\{G_i\}\longmapsto \varprojlim G_i$ induces an exact equivalence from $\Pro(\calF)$
to the category of profinite $k$-groups, and we henceforth identify these two categories.   

The functor  $\varpi_0\colon\calC\to\calF$ given by $G\mapsto G/G^0_{\red}$ is a left adjoint of the inclusion functor $\calF\hookrightarrow\calC$, hence $\varpi_0$ is right exact. It extends to a right exact functor on the associated pro-categories: 
\[
\varpi_0\colon\Pro(\calC)\longrightarrow\Pro(\calF)\quad\text{given by}\quad G=\{G_i\}\longmapsto G/G^0_{\red}\Def\{G_i/(G_i)^0_{\red}\}
\]
(cf. \cite[II, \S7]{Oort}). We may therefore consider the left-derived functors:
\[
\varpi_i\Def{\rm L}^i\varpi_0\colon\Pro(\calC)\longrightarrow\Pro(\calF),\quad i=0,1,2,\dots 
\]
For any object $G$ of $\Pro(\calC)$, we call $\varpi_1(G)$ the {\it Serre--Oort fundamental group} of $G$. 
According to \cite[II, \S14, Theorem (14.2)]{Oort}, we have $\varpi_i=0$ for $i\ge 2$. Hence, $\varpi_1$ is a left exact functor. Therefore, for any short exact sequence
\[
\xymatrix{
0\ar[r]& N\ar[r]& G\ar[r]&H\ar[r]& 0
}
\]
in $\Pro(\calC)$, we have a {\it homotopy exact sequence} in $\Pro(\calF)$:
\begin{equation}\label{eq:homotopy seq}
\xymatrix{
0\ar[r]&\varpi_1(N)\ar[r]&\varpi_1(G)\ar[r]&\varpi_1(H)\ar[r]&\varpi_0(N)\ar[r]&\varpi_0(G)\ar[r]&\varpi_0(H)\ar[r]& 0.
}\end{equation}

\begin{lem}\label{lem:varpi}
Let $G\in\Pro(\calF)$ and $H\in\Pro(\calC)$.
\begin{enumerate}
\item \label{lem:G-projective} $G$ is projective in $\Pro(\calF)$ if and only if it is projective
in $\Pro(\calC)$.

\item \label{lem:varpi-vanishing} We have $\varpi_1(G)=0$ and a natural isomorphism
$\varpi_1(H^0_{\red})\xrightarrow{\simeq}\varpi_1(H)$.
\end{enumerate}
\end{lem}

\begin{proof}
For \eqref{lem:G-projective}, we refer to \cite[II, \S7, Corollary (7.5)]{Oort}. To prove \eqref{lem:varpi-vanishing}, let $P_0\twoheadrightarrow G$ be a surjective morphism with $P_0$ projective in $\Pro(\calC)$. Then $\varpi_0(P_{0})$ is projective in $\Pro(\calF)$ by the universal property of $\varpi_0$. 
By replacing $P_0$ with $\varpi_0(P_0)$, we may assume that $P_{0}$ belongs to the category $\Pro(\calF)$. Then, as the kernel of the surjection $P_0\twoheadrightarrow G$ also belongs to $\Pro(\calF)$, repeating this construction inductively yields a projective resolution $P_{\bullet}$ of $G$ in $\Pro(\calF)$. 
By \eqref{lem:G-projective}, $P_{\bullet}$ is a projective resolution of $G$ in $\Pro(\calC)$. As $\varpi_0(P_\bullet)=P_{\bullet}$, we see that $\varpi_1(G)=0$. The second assertion then follows from the homotopy exact sequence \eqref{eq:homotopy seq} arising from the exact sequence $0\to H_{\red}^0\to H\to\varpi_0(H)\to 0$. 
\end{proof}

\begin{lem}(\cite[II, \S7]{Oort} and cf.\ \cite[\S6.2, Proposition 2]{Serre_GP})\label{lem:univ cov B}
For any object $H$ of $\Pro(\calC)$, 
there exists an extension in $\Pro(\calC)$
\begin{equation}
    \xymatrix{
        0 \ar[r] & \varpi_1(H) \ar[r] & \overline{H} \ar[r]^-{\overline{u}} & H^0_{\red} \ar[r] & 0
    }\label{eq:ext-univ}
\end{equation}
with the following universal property: any extension in $\Pro(\calC)$ of $H^0_{\red}$
by an object $G$ of $\Pro(\calF)$ arises by pushout of \eqref{eq:ext-univ}
along a unique homomorphism of profinite $k$-groups $\varpi_1(H)\rightarrow G$.  In other words, $\varpi_1(H)$ represents the functor 
$\Ext^1(H^0_{\red},-):\Pro(\calF)\rightarrow \Ab$.
\end{lem}

\begin{proof}
Let $P\twoheadrightarrow H^0_{\red}$ be a surjective morphism with $P$ projective. By applying the snake lemma to the commutative diagram with exact rows:
\[
\xymatrix{
0\ar[r]&\Ker(P_{\red}^0\to H^{0}_{\red})\ar[d]\ar[r]&P_{\red}^0\ar@{^{(}->}[d]\ar[r]&\im(P_{\red}^0\to H^{0}_{\red})\ar@{^{(}->}[d]\ar[r]&0\\
0\ar[r]&\Ker(P\twoheadrightarrow H^{0}_{\red})\ar[r]&P\ar[r]&H_{\red}^0\ar[r]&0,
}
\]
we get a surjection $\varpi_0(P)\twoheadrightarrow H_{\red}^0/\im(P_{\red}^0\to H^{0}_{\red})$. This together with the fact that $\varpi_0(H^0_{\red})=0$ implies that $H_{\red}^0/\im(P_{\red}^0\to H^{0}_{\red})=0$, hence the composition of maps $P^0_{\red}\hookrightarrow P\twoheadrightarrow H^0_{\red}$ is surjective. 

Set $R\Def\Ker(P^0_{\red}\to H^{0}_{\red})$. We put $\overline{H}\Def P^0_{\red}/R^0_{\red}$, so
by definition obtain an extension in $\Pro(\calC)$
\begin{equation}
 \xymatrix{
        0 \ar[r] & \varpi_0(R)=R/R^0_{\red}\ar[r] & \overline{H}=P^{0}_{\red}/R^0_{\red}\ar[r]^-{\overline{u}} & H^0_{\red} \ar[r] & 0 }.\label{eq:hbardefext}
\end{equation}
As $P/P^0_{\red}=\varpi_0(P)$ is in $\Pro(\calF)$, we have $\varpi_1(P/P^0_{\red})=0$ by Lemma \ref{lem:varpi}\eqref{lem:varpi-vanishing},
so   $\varpi_1(P^{0}_{\red})=\varpi_1(P)=0$ by the exact sequence \eqref{eq:homotopy seq} and the fact that $P$ is projective. 
Therefore, 
the homotopy exact sequence (\ref{eq:homotopy seq}) associated to the short exact sequence $0\to R^0_{\red}\to P^0_{\red}\to \overline{H}\to 0$ implies that $\varpi_1(\overline{H})=0$. We deduce that the homotopy exact sequence \eqref{eq:homotopy seq} associated to \eqref{eq:hbardefext} 
induces a natural isomorphism 
\begin{equation}
\varpi_1(H)=\varpi_1(H^0_{\red})\xrightarrow{~\simeq~}\varpi_0(R),
\label{eq:varpiHpi0R}
\end{equation}
where, as above, the first equality follows from  Lemma \ref{lem:varpi}\eqref{lem:varpi-vanishing}. This gives the extension \eqref{eq:ext-univ}.

If $G\in\Pro(\calF)$, applying $\Hom(-,G)$ to the tautological exact sequence
\[
\xymatrix{
0\ar[r]&P^0_{\red}\ar[r]&P\ar[r]&\varpi_0(P)\ar[r]&0
}
\]
yields the long exact sequence
\[
\xymatrix{
0\ar[r]&\Ext^1(\varpi_0(P),G)\ar[r]&\Ext^1(P,G)\ar[r]&\Ext^1(P^0_{\red},G)\ar[r]&\Ext^2(\varpi_0(P),G).
}
\]
As $P$ (respectively $\varpi_0(P)$) is projective in $\Pro(\calC)$ (respectively in $\Pro(\calF)$), it follows that 
\[
\Ext^1(P^{0}_{\red},G)=0
\]
for any $G\in\Pro(\calF)$. 
By construction, we have a surjection $\rho:P^0_{\red}\twoheadrightarrow \overline{H}$
with $\Ker(\rho)=R^0_{\red}$. 
It follows from this that $\Ext^1(\overline{H},G)=0$ for any $G\in\Pro(\calF)$.  
Applying $\Hom(-,G)$ to \eqref{eq:hbardefext}, we conclude that
the connecting homomorphism in the resulting long exact sequence
is a functorial isomorphism
\begin{equation}
\Hom(\varpi_1(H),G)\xrightarrow{~\simeq~}\Ext^1(H^0_{\red},G).
\label{eq:pi1 represents Ext}
\end{equation}
When 
$G=\varpi_1(H)$, the image of identity element ${\rm id}_{\varpi_1(H)}$ under \eqref{eq:pi1 represents Ext} is precisely the class of the 
extension \eqref{eq:ext-univ}; by Yoneda's lemma, 
\eqref{eq:hbardefext} therefore satisfies the stated universal property.
\end{proof}

\begin{definition}\label{def:univ-cover}
    For any object $H$ of $\Pro(\calC)$, we write $u:\overline{H}\to H$
    for the composition of the map $\overline{u}:\overline{H}\twoheadrightarrow H^0_{\red}$ in Lemma \ref{lem:univ cov B} with the inclusion $H^0_{\red}\hookrightarrow H$.
    We call $u:\overline{H}\to H$ the {\em universal cover}
    of $H$.  
\end{definition}

\begin{rem}\label{rem:univ-alt-character}
    It follows from the universal property of the extension 
    \eqref{eq:ext-univ} that the universal cover $u:\overline{H}\rightarrow H$
    is unique up to unique isomorphism.  Moreover, by construction we
    have
    \begin{enumerate}
    \item\label{eq:keru cokeru}$
        \Ker(u) \simeq \varpi_1(H),\ \Cok(u)\simeq \varpi_0(H)
    $
    \item \label{eq:varpi Hbar vaishing} $\varpi_i(\overline{H})=0$ for $i=0,1$.
    \end{enumerate}
    In fact, these properties also characterize the universal cover
    up to unique isomorphism.  
    Indeed,
    any pair $(\overline{H},u)$ satisfying \eqref{eq:keru cokeru}--\eqref{eq:varpi Hbar vaishing} yields an extension as in \eqref{eq:hbardefext}.
    As $\varpi_i(\overline{H})=0$ for $i=0,1$, we
    may choose a surjection
    $\rho:P\twoheadrightarrow \overline{H}$ with $P$ projective and both $P$ and $\Ker(\rho)$ connected and reduced.  As in the proof
    of Lemma \ref{lem:univ cov B}, one deduces that $\Ext^1(\overline{H},G)=0$ for any profinite $G$, and hence that the extension determined 
    by $(\overline{H},u)$ satisfies the universal property
    of Lemma \ref{lem:univ cov B}.
\end{rem}

\begin{rem}\label{rem:varpi}
When $k$ is an algebraically closed field of characteristic $0$, the isomorphism (\ref{eq:pi1 represents Ext}) is due to \cite[\S5.4, Corollaire]{Serre_GP}. In \cite[Corollary 2.7]{Brion}, Brion establishes an analogue of this in a more general setting.
\end{rem}

We now prove our main result:

\begin{thm}\label{thm:ab Nori}
There exists an isomorphism of profinite $k$-groups, 
\[
\pi^{\ab}(U)\xrightarrow{~\simeq~}\varprojlim_{\frkm}\varpi_1(J_{X,\frkm}),
\]
where $\frkm$ runs over the set of moduli supported on $S$. 
\end{thm}

\begin{proof}
For any modulus $\frkm$ supported on $S$, pullback
of the universal cover $\overline{J}_{X,\frkm}\to J_{X,\frkm}$ 
along the map $\varphi_{\frkm}\colon U\to J_{X,\frkm}$
gives a $\varpi_1(J_{X,\frkm})$-torsor $\varphi_{\frkm}^*\overline{J}_{X,\frkm}\to U$. 
For variable $\frkm$ supported on $S$, these give rise
to a homomorphism
\[
\pi^{\ab}(U)\longrightarrow\varprojlim_{\frkm}\varpi_1(J_{X,\frkm}).
\]
that we claim is an isomorphism. It suffices to show that for any finite $k$-group $G$, the induced map
\[
\varinjlim_{\frkm}\Hom(\varpi_1(J_{X,\frkm}),G)=\Hom\left(\varprojlim_{\frkm}\varpi_1(J_{X,\frkm}),G\right)\longrightarrow\Hom(\pi^{\ab}(U),G)
\]
is bijective, where the first identification is a consequence of \cite[Proposition 04AK]{SP}. By Theorem \ref{thm:CFT}, it suffices to show that the connecting homomorphism
\[
\Hom(\varpi_1(J_{X,\frkm}),G)\longrightarrow\Ext^1(J_{X,\frkm},G)
\]
is bijective for any $\frkm$. This is valid thanks to the isomorphism (\ref{eq:pi1 represents Ext}).  
\end{proof}

As a special case of Theorem \ref{thm:ab Nori}, we deduce the following well-known
description of the abelianized Nori fundamental group scheme of a projective curve:

\begin{cor}(cf.\ \cite[Remark 2]{Nori83}, \cite[Corollary 3.8]{Antei}, \cite[Corollary 7.2]{Langer})\label{cor:ab pi X}
Let $X$ be a smooth projective curve over $k$. The Albanese map $\varphi_0\colon X\to J_{X}$ associated to any fixed choice
of base point induces a natural isomorphism of profinite $k$-groups
\[
\varphi_{0,*}\colon\pi^{\ab}(X)\xrightarrow{~\simeq~}\varprojlim_{n}J_X[n].
\]
\end{cor}

To prove this, we need: 

\begin{lem}\label{lem:cor:ab pi X}
For any abelian variety $A$ over $k$, we have a natural identification
\[
\varpi_1(A)\xrightarrow{~\simeq~}\varprojlim_nA[n].
\]
\end{lem}

\begin{proof}
For any subgroup $N$ of $\varpi_1(A)$ with $\varpi_1(A)/N$ a finite $k$-group, we have $n(\varpi_1(A)/N)=0$, and hence $n\varpi_1(A)\subset N$ for all $n>0$ sufficiently large with respect to the divisibility ordering on the positive integers 
\cite[VIIA, Proposition 8.5]{SGA3i}. 
It follows that the canonical map
\begin{equation}
\varpi_1(A)\xrightarrow{~\simeq~}\varprojlim_n\varpi_1(A)/n\varpi_1(A)
\label{eqn:varpi-Jac-A}
\end{equation}
is an isomorphism of profinite $k$-groups.
On the other hand, for any integer $n>0$, the exact sequence
\[
\xymatrix{
0\ar[r]& A[n]\ar[r]& A\ar[r]& A\ar[r]& 0
}\]
gives rise to the long exact homotopy sequence
\[
\xymatrix{
0\ar[r]& \varpi_1(A[n])\ar[r]& \varpi_1(A)\ar[r]^{n}& \varpi_1(A)\ar[r]&\varpi_0(A[n])\ar[r]& 0.
}
\]
As $A[n]$ is a finite $k$-group, we have $\varpi_0(A[n])=A[n]$, and $\varpi_1(A[n])=0$ thanks to Lemma \ref{lem:varpi}(2). We obtain an exact sequence
\begin{equation}\label{eqn:varpi-Jac-B}
\xymatrix{
0\ar[r]&\varpi_1(A)\ar[r]^{n}&\varpi_1(A)\ar[r]& A[n]\ar[r]& 0.
}
\end{equation}
Combining \eqref{eqn:varpi-Jac-A} with \eqref{eqn:varpi-Jac-B}
yields the desired identification. 
\end{proof}

\begin{proof}[Proof of Corollary \ref{cor:ab pi X}]
Specializing Theorem \ref{thm:ab Nori} to the case $U=X$ and $S=\emptyset$,
we see that $\varphi_0$ induces an isomorphism of profinite $k$-groups $\varphi_{0,*}\colon\pi^{\ab}(X)\xrightarrow{\simeq}\varpi_1(J_X)$. Therefore, the corollary follows from Lemma \ref{lem:cor:ab pi X}.
\end{proof}

\begin{cor}\label{cor:mult}
The multiplicative part $\pi^{\ab}(U)^{\rm m}$ of the abelian fundamental group scheme can be described as follows:
\[
\pi^{\ab}(U)^{\rm m}\simeq\Diag\left(J_X(k)_{\rm tor}\times\Ker\left(\Sigma\colon\prod_{x\in S}\Q/\Z\to\Q/\Z\right)\right).
\]
\end{cor}

\begin{proof}
As $k$ is algebraically closed, a $k$-group is of multiplicative type if and only if it is diagonalizable, i.e., $G\simeq\Diag(M)$ for some abelian group $M$ \cite[\S 7.3, Corollary]{Wat}, in which case we have $M\simeq\Hom_{\text{$k$-grp}}(G,\G_m)$, and 
the multiplicative part $\pi^{\ab}(U)^{\rm m}$ is defined to be
\[
\pi^{\ab}(U)^{\rm m}\Def\Diag\left(\Hom(\pi^{\ab}(U),\G_m)\right).
\]
Therefore, it suffices to show that
\[
\Hom(\pi^{\ab}(U),\G_m)\simeq J_X(k)_{\rm tor}\times\Ker\left(\Sigma\colon\prod_{x\in S}\Q/\Z\to\Q/\Z\right).
\]
By Theorem \ref{thm:ab Nori}, we have
\begin{align*}
\Hom(\pi^{\ab}(U),\G_m)&=\varinjlim_{\frkm}\Hom(\varpi_1(J_{X,\frkm}),\G_m)=\Hom(\varpi_1(J_{X,\frkm_{\red}}),\G_m)\\
&=\varinjlim_{n}\Hom(\varpi_1(J_{X,\frkm_{\red}}),\mu_n)=\varinjlim_n\Ext^1(J_{X,\frkm_{\red}},\mu_n)
\end{align*}
By Remark~\ref{rem:L_m}, the generalized Jacobian $J_{X,\frkm}$ is naturally an extension of $J_X$ by $(\prod_{x\in S}\G_m)/\G_m$;
applying $\varinjlim_n\Ext^1(\cdot,\mu_n)$ to the corresponding short exact sequence gives 
\begin{equation*} 
\xymatrix{
	0 \ar[r] & \varinjlim_n \Ext^1(J_X,\mu_n) \ar[r] & \Hom(\pi^{\ab}(U),\G_m) \ar[r] & \varinjlim_{n} \Ext^1\left(\left(\prod_{x\in S}\G_m\right)/\G_m,\mu_n\right) \ar[r] & 0
}.
\end{equation*}
Via the natural Kummer theory identifications
\[
\Ext^1(J_X,\mu_n)\simeq J_X[n](k)\qquad\text{and}\qquad \Ext^1(\G_m,\mu_n) \simeq \Z/n\Z,	
\]
the outer terms of the sequence are identified with $J_X(k)_{\mathrm{tor}}$ and
$\Ker(\Sigma:\prod_{x\in S}\Q/\Z\to\Q/\Z)$, respectively.
Since $J_X(k)_{\mathrm{tor}}$ is divisible, the sequence splits
noncanonically, yielding the asserted description.
\end{proof}

\begin{cor}\label{cor:pro-p et}
The maximal pro-$p$ quotient $\pi^{\ab}(U)(k)^{(p)}$ of the profinite group $\pi^{\ab}(U)(k)$ admits the following description as a topological group:
\[
\pi^{\ab}(U)(k)^{(p)}\simeq T_pJ_X(k)\times 
\prod_{x\in S}\prod_{i\ge 1,\,p\nmid i}\Hom\left(CW(k),\Q_p/\Z_p\right),
\]
in which $T_pJ_X(k):=\varprojlim_nJ_X[p^n](k)$ and $CW(k):=\varinjlim_s(W_s(k),V)$ is the group of Witt-covectors
equipped with the discrete topology, whose
character group $\Hom(CW(k),\Q_p/\Z_p)$ carries the topology of pointwise
convergence.
\end{cor}

Before giving the proof, we recall some preliminary facts from \cite[\S8]{Serre_GP}. For ease of notation, set $C:=\Q_p/\Z_p$ and $P_s\:=\varpi_1(W_s)(k)$.
Serre's self-duality \cite[\S8.4, Proposition 4]{Serre_GP}, together with (\ref{eq:pi1 represents Ext}), gives isomorphisms
\[
\theta_s:W_s(k)\xrightarrow{~\simeq~}
\Hom_{\rm cont}(P_s,C)
\quad\text{satisfying the compatibility}\quad R^*\theta_s=\theta_{s+1} F V,
\]
where $R:W_{s+1}\to W_s$ is truncation and $R^*$ 
is the pullback map on $C$-valued characters induced by
$\varpi_1(R)(k):P_{s+1}\to P_s$ (see \cite[\S8.4, Lemma 2]{Serre_GP}).
To obtain transition maps involving $V$ alone, we replace
$\theta_s$ by
\begin{equation}\label{eq:serre theta}
\widetilde\theta_s:=\theta_s  F^s,
\end{equation}
which is again an isomorphism as $F$ is bijective on $W_s(k)$ due to perfectness of $k$.
Since $F$ commutes with $V$, we then have the compatibility relation $R^*\widetilde\theta_s = \widetilde\theta_{s+1}V$,
and the dual identifications
\[
P_s\xrightarrow{~\simeq~}\Hom(W_s(k),C),
\qquad
g\longmapsto\bigl(a\mapsto\widetilde\theta_s(a)(g)\bigr),
\]
intertwine $\varpi_1(R)(k):P_{s+1}\to P_s$ with the map
\begin{equation}\label{eq:serre theta V}
\Hom(W_{s+1}(k),C)\longrightarrow\Hom(W_s(k),C),
\qquad \chi\longmapsto\chi\circ V.
\end{equation}
We will use these facts and the above notation in the proof below.

\begin{proof}[Proof of Corollary \ref{cor:pro-p et}]
For $S=\emptyset$, the claimed identification follows from Corollary \ref{cor:ab pi X}, so we henceforth assume $S\ne\emptyset$.
By Theorem \ref{thm:ab Nori}, we have
\begin{align*}
\pi^{\ab}(U)(k)^{(p)}&=\varprojlim_{\frkm}\varpi_1(J_{X,\frkm})(k)^{(p)}\simeq\varprojlim_{\frkm,n}\varpi_1(J_{X,\frkm})(k)/p^n\varpi_1(J_{X,\frkm})(k).
\end{align*}
For any modulus $\frkm$ supported on $S$ and any integer $n>0$, consider the commutative diagram with exact rows: 
\[
\xymatrix{
0\ar[r]&\varpi_1(L_{\frkm})\ar[r]\ar[d]^{p^n}&\varpi_1(J_{X,\frkm})\ar[r]\ar[d]^{p^n}&\varpi_1(J_X)\ar[r]\ar@{^{(}->}[d]^{p^n}&0\\
0\ar[r]&\varpi_1(L_{\frkm})\ar[r]&\varpi_1(J_{X,\frkm})\ar[r]&\varpi_1(J_X)\ar[r]&0
}
\]
By the snake lemma, together with Lemma \ref{lem:cor:ab pi X}, we have the exact sequence of $\Z/p^n\Z$-modules
\[
\xymatrix{
0\ar[r]&\varpi_1(L_{\frkm})(k)/p^n\varpi_1(L_{\frkm})(k)\ar[r]&\varpi_1(J_{X,\frkm})(k)/p^n\varpi_1(J_{X,\frkm})(k)\ar[r]& J_X[p^n](k)\ar[r]& 0.
}
\]
As the inverse system $\{\varpi_1(L_{\frkm})(k)/p^n\varpi_1(L_{\frkm})(k)\}_{n\ge 1}$ is Mittag--Leffler, by \cite[Lemma 0598]{SP}, we obtain a short exact sequence of $\Z_p$-modules
\[
\xymatrix{
0\ar[r]&\varprojlim_n\varpi_1(L_{\frkm})(k)/p^n\varpi_1(L_{\frkm})(k)\ar[r]&\varprojlim_{n}\varpi_1(J_{X,\frkm})(k)/p^n\varpi_1(J_{X,\frkm})(k)\ar[r]& T_pJ_X(k)\ar[r]& 0.
}
\]
Taking the inverse limit over $\frkm$, which is exact for
compact abelian groups, gives
\[
0\longrightarrow
\varprojlim_{\frkm,n}
\varpi_1(L_{\frkm})(k)/p^n\varpi_1(L_{\frkm})(k)
\longrightarrow \pi^{\ab}(U)(k)^{(p)}
\longrightarrow T_pJ_X(k)\longrightarrow0.
\]
As $T_pJ_X(k)$ is a finite free
$\Z_p$-module, this sequence (noncanonically) splits, giving an isomorphism
\begin{equation}\label{eq:pro-p et}
\pi^{\ab}(U)(k)^{(p)}
\simeq T_pJ_X(k)\times
\varprojlim_{\frkm,n}\varpi_1(L_{\frkm})(k)/p^n\varpi_1(L_{\frkm})(k).
\end{equation}

It remains to describe
$\varprojlim_n\varpi_1(L_{\frkm})(k)/
p^n\varpi_1(L_{\frkm})(k)$,
and we may restrict to moduli $\frkm$ with ${\rm Supp}(\frkm)=S$
as these form a cofinal subsystem. 
Let $T_{\frkm}$ (respectively $U_{\frkm}$) be the multiplicative part (respectively the unipotent part) of $L_{\frkm}$. 
Writing $\frkm=\sum_{x\in S}n_x x$, the description of $L_{\frkm}$ recalled in Remark \ref{rem:L_m} gives
$
U_{\frkm}=\prod_{x\in S}V_{(n_x)}.
$
Now for any positive integer $m$, the group  $V_{(m)}$ is a finite product of finite-length Witt groups by Remark \ref{rem:L_m};
in particular, $U_{\frkm}$ is connected and reduced.  As the multiplication by $p^n$ map on $T_{\frkm}$ is finite surjective,
we deduce the identifications
\[
\varpi_1(L_{\frkm})(k)/p^n\varpi_1(L_{\frkm})(k)=\varpi_1(L_{\frkm}/p^nL_{\frkm})(k)=\varpi_1(U_{\frkm}/p^nU_{\frkm})(k).
\]
On the other hand, the explicit description of $V_{(n_x)}$ in (\ref{eq:local grp}) gives
\[
U_{\frkm}/p^nU_{\frkm}=\prod_{x\in S}V_{(n_x)}/p^nV_{(n_x)}=\prod_{x\in S}\prod_{1\le i<n_x,p\nmid i}W_{\min\{n,r_i(n_x)\}},
\]
where $r_i(n_x):=\min\{r>0\,|\,ip^r\ge n_x\}$, hence
\begin{align*}
\varpi_1(L_{\frkm})(k)/p^n\varpi_1(L_{\frkm})(k)&=\prod_{x\in S}\prod_{1\le i<n_x,p\nmid i}P_{\min\{n,r_i(n_x)\}}\xrightarrow{\simeq}\prod_{x\in S}\prod_{1\le i<n_x,p\nmid i}\Hom(W_{\min\{n,r_i(n_x)\}}(k),C),
\end{align*}
where the second identification is induced by $\widetilde{\theta}_{\min\{n,r_i(n_x)\}}$ given by (\ref{eq:serre theta}). 
Under these identifications, the
natural transition map
\[
\varpi_1(L_{\frkm})(k)/p^{n+1}\varpi_1(L_{\frkm})(k)
\longrightarrow
\varpi_1(L_{\frkm})(k)/p^n\varpi_1(L_{\frkm})(k)
\]
corresponds to the map induced by the morphism
\[
\Hom(W_{\min\{n+1,r_i(n_x)\}}(k),C)\longrightarrow \Hom(W_{\min\{n,r_i(n_x)\}}(k),C)\]
which is the identity if $r_i(n_x)\le n$ and restriction along
$V\colon W_n\hookrightarrow W_{n+1}$ if $r_i(n_x)>n$ (see (\ref{eq:serre theta V})). 
We deduce
\[
\varprojlim_n\varpi_1(L_{\frkm})(k)/p^n\varpi_1(L_{\frkm})(k)\simeq \prod_{x\in S}\prod_{1\le i<n_x,p\nmid i}\Hom\left(W_{r_i(n_x)}(k),C\right).
\]
Combining this with (\ref{eq:pro-p et}), we obtain
\begin{align*}
\pi^{\ab}(U)(k)^{(p)}&\simeq
T_pJ_X(k)\times\varprojlim_{\frkm}
\prod_{x\in S}\prod_{1\le i<n_x,p\nmid i}\Hom\left(W_{r_i(n_x)}(k),C\right)\\
&=T_pJ_X(k)\times
\prod_{x\in S}\prod_{i\ge 1\,,p\nmid i}\Hom\left(\varinjlim_{n_x}W_{r_i(n_x)}(k),C\right)\\
&=T_pJ_X(k)\times\prod_{x\in S}\prod_{i\ge 1,\,p\nmid i}\Hom(CW(k),C).
\end{align*}
Here, the last identification follows from (\ref{eq:serre theta V}) together with the fact that,
for each fixed $i$ prime to $p$, the length $r_i(n_x)$
tends to infinity as $n_x\to\infty$. 
This completes the proof. 
\end{proof}

\begin{rem}
Note that $\pi^{\ab}(U)(k)$ is simply the maximal abelian quotient of the \'etale fundamental group $\pi_1^{\et}(U)$ of the curve $U$. 
With this in mind, the explicit description given by Corollary \ref{cor:pro-p et} can be deduced via Pontryagin duality from
Katz's local-to-global sequence for \'etale characters
\cite[\S2.2, (2.2.2.1)]{Katz86},
using \cite[Theorem~4.9(i), Remark~1.1]{CP09} to express each local
character group as a direct sum of Witt covectors.
On the other hand, in \cite[Theorem 3.1, Corollary 3.5]{Kumar08}, Kumar provides an explicit description of the maximal pro-$p$ abelian quotient of the \'etale fundamental group $\pi_1^{\et}(U)$ of an affine curve $U$ in terms of its coordinate ring
$A:=\Gamma(U,\scrO_U)$.  By Artin--Schreier--Witt theory,  one has the identification
\[
H^1_{\et}(U,\Q_p/\Z_p)
\simeq\varinjlim_n W_n(A)/(F-1)W_n(A),
\]
with transition maps induced by $V$, and Kumar's description is the Pontryagin dual of this.
\end{rem}

For any $k$-group $G$ and any integer $m$, we write
$G^{(m)}$ for the pullback of $G$ along the automorphism
$F^m\colon k \rightarrow k$ of $k$.  By definition,
we have $(G^{(m)})^{(1)}=G^{(m+1)}$, and for any 
nonnegative integer $n$, the $n$-th iterate of the relative
Frobenius map defines a morphism $F^n: G^{(m)}\rightarrow G^{(m+n)}$.
We denote by $G[F^n]$ the kernel of $F^n:G\rightarrow G^{(n)}$.
If $P_{\bullet}$ is a projective resolution of $G$, then the relative Frobenius morphism $F$ on $G$ is extended by the relative Frobenius morphism
on each term of the resolution. Noting that the Frobenius twist of the resolution $P_{\bullet}^{(1)}$ gives a projective resolution of $G^{(1)}$, by applying $\varpi_0(-)$
to the resolution and passing to homology, we deduce that
the morphism $\varpi_1(F)$ coincides with the relative Frobenius
morphism of $\varpi_1(G)$.

\begin{lem}\label{lem:conn part}
For any connected and reduced $k$-group $G$, we have
\[
\varpi_1(G)^0\simeq\varprojlim_nG^{(-n)}[F^n],
\]
with transition maps induced by the relative Frobenius morphism $F\colon G^{(-n-1)}\to G^{(-n)}$. 
\end{lem}

\begin{proof}
For any $k$-subgroup $N$ of $\varpi_1(G)^0$ with $\varpi_1(G)^0/N$ finite, as $\varpi_1(G)^0/N$ is finite 
 and connected, we have $F^n(\varpi_1(G)^{0,(-n)}/N^{(-n)})=0$, hence $F^n(\varpi_1(G)^{0,(-n)})\subset N$ for sufficiently large $n>0$.
This implies that we have a natural isomorphism 
\[
\varpi_1(G)^0\xrightarrow{~\simeq~}\varprojlim_n\varpi_1(G)^0/F^n(\varpi_1(G)^{0,(-n)}).
\]
Per the discussion above, for any positive integer $n>0$, the commutative diagram with exact rows
\[
\xymatrix{
0\ar[r]& G^{(-n-1)}[F^{n+1}]\ar[r]\ar[d]_F& G^{(-n-1)}\ar[r]_{F^{n+1}}\ar[d]_F& G\ar[r]\ar@{=}[d]& 0\\
0\ar[r]& G^{(-n)}[F^n]\ar[r]& G^{(-n)}\ar[r]_{F^n}& G\ar[r]& 0
}
\]
gives rise to the commutative diagram with exact rows
\[
\xymatrix{
0\ar[r]& \varpi_1(G)^{(-n-1)}\ar[r]_{F^{n+1}}\ar[d]_F& \varpi_1(G)\ar[r]\ar@{=}[d]& G^{(-n-1)}[F^{n+1}]\ar[r]\ar[d]_{F}& 0\\
0\ar[r]& \varpi_1(G)^{(-n)}\ar[r]_{F^n}& \varpi_1(G)\ar[r]& G^{(-n)}[F^n]\ar[r]& 0
}
\]
from which we deduce a natural isomorphism  
\[
\varpi_1(G)^0\xrightarrow{~\simeq~}\varprojlim_n\varpi_1(G)^0/F^n(\varpi_1(G)^{0,(-n)})\xrightarrow{~\simeq~}\varprojlim_{n}G^{(-n)}[F^n].
\qedhere
\]
\end{proof}

\begin{cor}\label{cor:conn part}
The connected component of the identity $\pi^{\ab}(U)^0$ of the abelian fundamental group scheme can be described as follows:
\[
\pi^{\ab}(U)^0\xrightarrow{~\simeq~}
\varprojlim_{\frkm}\varprojlim_{n}J^{(-n)}_{X,\frkm}[F^n]=
\varprojlim_{\frkm,n}J_{X,\frkm}^{(-n)}[F^n],
\]
where the transition maps are induced by the relative Frobenius morphism $F\colon J_{X,\frkm}^{(-n-1)}\to J_{X,\frkm}^{(-n)}$. 
\end{cor}

\begin{proof}
This follows immediately from Theorem \ref{thm:ab Nori}
and Lemma \ref{lem:conn part}.
\end{proof}

\section{The unipotent fundamental group schemes for singular curves}\label{sec:singular curve}

It does not seem that there exists a relevant description of the unipotent part of the abelian fundamental group scheme $\pi^{\ab}(U)$ purely in terms of the generalized Jacobians $J_{X,\frkm}$ unless $U=X$. In fact, the following result due to Nori suggests that each $J_{X,\frkm}$ is directly related to the abelian unipotent fundamental group scheme of the corresponding singular curve $X_{\frkm}$ rather than the affine curve $U$. The curve $X_{\frkm}$ is defined in \cite[IV,\S1]{Serre}; it is a projective and integral curve over $k$ with normalization $X\rightarrow X_{\frkm}$ that is an isomorphism
over $X\setminus S$, 
and it satisfies $\Pic_{X_{\frkm}}^0=J_{X,\frkm}$.

\begin{thm}(\cite[IV, Proposition 6]{Nori82})\label{thm:sing curve}
The unipotent part of $\pi^{\ab}(X_{\frkm})$ can be described as follows:
\[
\pi^{\ab}(X_{\frkm})^{\uni}=\varprojlim_{G=G^0\hookrightarrow J_{X,\frkm}}G^D=\varprojlim_{n}J_{X,\frkm}[F^n]^D,
\]
where $G^D$ means the Cartier dual of $G$, and the limit runs over all finite local subgroup schemes of the generalized Jacobian $J_{X,\frkm}$. 
\end{thm}

\begin{rem}
Note that the map $\pi^{\ab}(U)\to\pi^{\ab}(X_{\frkm})$ induced by the composition $U\hookrightarrow X\to X_{\frkm}$ is far from being injective or surjective. For example, when $X=\bbP^1_k$, the map $\pi^{\ab}(U)\to\pi^{\ab}(X_{\frkm})$ factors through the trivial group scheme  $\pi^{\ab}(\bbP^1_k)=1$. 
\end{rem}

We conclude the present paper with a proof of the following version of Das' description of the \'etale fundamental group of a singular curve \cite[Theorem 3.3]{Das}: 

\begin{thm}\label{thm:uni ab}
Let $\frkm=\sum_{x\in S}n_xx$ be a nonzero modulus with support $T:={\rm Supp}(\frkm)\subseteq S$. For each point $x\in T$ and any $1\le i\le n_x-1$ with $p\nmid i$, we denote by $r_{x,i}$ the smallest integer $r>0$ satisfying $p^r\ge n_x/i$. Then $\pi^{\ab}(X_{\frkm})^{\uni}$ fits into the short exact sequence 
\[
\xymatrix{
0\ar[r]&\pi^{\ab}(X)^{\uni}\ar[r]&\pi^{\ab}(X_{\frkm})^{\uni}\ar[r]&\Z_p^{\oplus \#T-1}\times\prod_{x\in T}\prod_{1\le i \le n_x-1,\,p\nmid i} W[F^{r_{x,i}}]\ar[r]& 0,
}
\]
in which $W\Def\varprojlim_n W_n$. 
\end{thm}

\begin{proof}
Taking Frobenius kernels on the exact sequence
$0\to L_{\frkm}\to J_{X,\frkm}\to J_X\to 0$,
yields exact sequences 
\[
\xymatrix{
0\ar[r]& L_{\frkm}[F^n]\ar[r]& J_{X,\frkm}[F^n]\ar[r]& J_X[F^n]\ar[r]& 0
}\]
for all $n\ge 1$, which form an inductive system of short exact sequences.  
Dualizing, we obtain the exact sequences
\[
\xymatrix{
0\ar[r]& J_X[F^n]^D\ar[r]& J_{X,\frkm}[F^n]^D\ar[r]& L_{\frkm}[F^n]^D\ar[r]& 0
}
\]
for all $n\ge 1$, which form a projective system of short exact sequences of finite $k$-groups. 
Taking the projective limit of these, 
we obtain an exact sequence of profinite $k$-groups
\[
\xymatrix{
0\ar[r]& \varprojlim_nJ_X[F^n]^D\ar[r]&\varprojlim_n J_{X,\frkm}[F^n]^D\ar[r]&\varprojlim_nL_{\frkm}[F^n]^D\ar[r]& 0\\
&\pi^{\ab}(X)^{\uni}\ar[u]_{\simeq}&\pi^{\ab}(X_{\frkm})^{\uni}\ar[u]_{\simeq}&&
}
\]
where the vertical isomorphisms are provided by Theorem \ref{thm:sing curve}. Observe that, for each $n\ge 1$, the map $J_{X,\frkm}[F^n]^D\to L_{\frkm}[F^n]^D$ is a faithfully flat morphism between finite $k$-groups. Since faithful flatness for a morphism between affine $k$-group schemes is equivalent to injectivity of the corresponding Hopf algebra homomorphism, and filtered colimits of Hopf algebra homomorphisms preserve injectivity, the induced morphism of profinite $k$-groups $\varprojlim_n J_{X,\frkm}[F^n]^D\to\varprojlim_n L_{\frkm}[F^n]^D$ is faithfully flat. Hence the above sequence is right exact.

It therefore suffices to show that
\[
\varprojlim_nL_{\frkm}[F^n]^D\simeq\Z_p^{\oplus\#T-1}\times\prod_{x\in T}\prod_{1\le i \le n_x-1,\,p\nmid i} W[F^{r_{x,i}}].
\]
To that end, it suffices to show that 
\begin{equation}\label{eq:uni sing}
L_{\frkm}[F^n]^D\simeq(\Z/p^n\Z)^{\oplus\#T-1}\times\prod_{x\in T}\prod_{1\le i \le n_x-1,\,p\nmid i} W_n[F^{r_{x,i}}]
\end{equation}
compatibly in each $n\ge 1$. Note that $\Z/p^n\Z=\mu_{p^n}^D$ and $W_n[F^{r_{x,i}}]=W_{r_{x,i}}[F^n]^D$ (see \cite[II,\S10]{Oort} for the latter identification; this can be proved via Dieudonn\'e theory by observing that for the contravariant Dieudonn\'e module functor $\mathbb{D}$, we have a natural identification $\mathbb{D}(W_n[F^m]) = W(k)[F,V]/(F^m,V^n)$). 
Thus, the equality (\ref{eq:uni sing}) is immediate from the description of $L_{\frkm}$ given by \eqref{eq:local grp} and Remark \ref{rem:L_m}.
\end{proof}



\printbibliography

\end{document}